\documentclass[preprint,12pt]{elsarticle}
\usepackage{eurosym}
\usepackage{amssymb}
\usepackage{array}
\usepackage{booktabs}
\usepackage{amsmath}
\usepackage{lineno}
\usepackage{latexsym,amscd,amsthm,amsfonts,amstext}
\usepackage[mathscr]{eucal}
\usepackage{graphicx}
\usepackage{subfig}
\usepackage{float}
\usepackage{color}
\usepackage{hyperref}
\usepackage[utf8]{inputenc}
\usepackage[english]{babel}

\newtheorem{assumption}{Assumption}[section]
\newtheorem{theorem}{Theorem}[section]

\newtheorem{remark}{Remark}[section]
\numberwithin{equation}{section}

\journal{Computers \& Mathematics with Applications}

\begin{document}

\begin{frontmatter}



\title{Toward Practical Forecasts of Public Sentiments via Convexification for Mean Field Games: Evidence from Real World COVID-19 Discussion Data}



\author[1]{Shi Chen}
\author[2]{Michael V. Klibanov\texorpdfstring{\corref{cor}}{}} \ead{mklibanv@charlotte.edu}
\author[2]{Kevin McGoff}
\author[3]{Trung Truong}
\author[4]{Wangjiaxuan Xin}
\author[1]{Shuhua Yin}

\affiliation[1]{organization={Department of Epidemiology and Community Health, University of North Carolina at Charlotte},
            city={Charlotte},
            postcode={28223}, 
            state={NC},
            country={USA}}

\affiliation[2]{organization={Department of Mathematics and Statistics, University of North Carolina at Charlotte},
            city={Charlotte},
            postcode={28223}, 
            state={NC},
            country={USA}}

\affiliation[3]{organization={Department of Mathematics and Physics, Marshall University},
            city={Huntington},
            postcode={25755}, 
            state={WV},
            country={USA}}
\affiliation[4]{organization={Department of Software and Information Systems, University of North Carolina at Charlotte},
            city={Charlotte},
            postcode={28223}, 
            state={NC},
            country={USA}}

\cortext[cor]{corresponding author}

\begin{abstract}
We apply a convexification-based numerical method to forecast public sentiment dynamics using Mean Field Games (MFGs). The theoretical foundation for the convexification approach, established in our prior work, guarantees global convergence to the unique solution to the MFG system. The present work demonstrates the practical potential of this framework using real-world sentiment data extracted from social media public discussion during the COVID-19 pandemic. The results show that the MFG model with appropriate parameters and convexification
yields sentiment density predictions that align closely with observed data and satisfy the governing equations. While current parameter selection relies on manual calibration, our findings establish the first proof-of-concept evidence that MFG models can capture complex temporal patterns in public sentiment, laying the groundwork for future work on systematic parameter identification methods, i.e. solutions of coefficient inverse problems for the MFG system.
\end{abstract}




\end{frontmatter}



\section{Introduction}

Mean Field Games (MFGs) theory, introduced independently by Lasry and Lions 
\cite{LasryLions2007} and by Huang, Caines, and Malham\'e \cite%
{HuangCainesMalhame2006}, provides a powerful framework for modeling
large-scale interacting agent systems. The theory has found applications in
diverse domains, e.g., economics \cite{Achdou,Achdou2}, finance \cite%
{LasryLions2007,Cardaliaguet,Firoozi}, crowd dynamics \cite{Dogbe,Lachapelle}%
, and opinion formation \cite{Bauso,Banez, Festa, Gao, Stella}.

Forecasting public sentiment, especially its dynamic shift, is a challenging problem with significant implications for public health policy, crisis management, and social
planning during major emergencies, including the COVID-19 pandemic. Proactive forecast of public sentiment provides an effective way for public health practitioners to evaluate potential policy adoptions by the public and to further estimate epidemic burdens in the society. Traditional statistical and machine learning approaches often treat sentiment prediction as a technical data-fitting exercise without incorporating the underlying dynamical structure of public sentiment formation and shift. In
contrast, MFG-based models respect fundamental principles of agents
interactions, making them generalizable to unseen data. Prior work has
applied MFGs to opinion dynamics and sentiment analysis \cite{Banez, Festa,
Gao, Stella}, but none of these works provides a solver with rigorously justified global
convergence for sentiment forecasting. Our convexification approach in \cite%
{Kliba2025} offers such a guarantee by transforming the original problem
into a convex optimization problem using the Carleman Weight Functions
(CWFs), thereby avoiding the challenges of local minima that often plague
traditional methods.

The convexification method was first introduced in \cite{Klib95,Klib97} for
two coefficient inverse problems (CIPs) for hyperbolic partial differential
equations (PDEs). The main purpose of the convexification is to handle the well
known phenomenon of multiple local minima and ravines of conventional least
squares cost functionals for CIPs, which arises since these functionals are not convex in general. The
convexification method works for many CIPs. While the works \cite{Klib95,Klib97}
are purely theoretical, more recent publications contain both theory
and numerical studies, see \cite{KL,Kepid} for some examples. In addition, convexification is applicable to some CIPs for MFGs \cite{MFGbook,MFGCAMWA,MFG8IPI}.

The convexification method is a numerical version of the theoretical
publication \cite{BukhKlib}. In \cite{BukhKlib}, the apparatus of Carleman
estimates was introduced in the field of Inverse Problems. The
convexification constructs a weighted Tikhonov-like functional for a CIP.
The weight is the CWF. This function is used as the
weight in the Carleman estimate for the corresponding PDE operator. The key
result states that this functional is strongly convex on an appropriate
convex bounded set. The diameter $d>0$ of this set is an arbitrary one.
Also, that functional has unique minimizer on that set. Carleman estimates
were introduced in MFGs in \cite{MFG1}. There were a number of works since
then, which have developed this technique further for MFGs, see, e.g. \cite{KL,MFGCAMWA,MFG8IPI,MFGSAM,MFG2}. In addition to CIPs, it was established in our previous work \cite{Kliba2025}
that the convexificaton concept can be applied to the problem of forecasting
public sentiments via MFGs. Specifically, convexification has been justified theoretically and validated on simulated data in \cite{Kliba2025}.

In this paper, we present the first demonstration of the practical viability of the convexification
framework for public sentiment forecasting using real-world data. Our dataset
consists of dynamic, daily discussions of the general public regarding COVID-19 on the social media platform X (formerly Twitter\footnote{%
\url{https://www.X.com/}}) over a two-year period during the pandemic (March 2020--April 2022).
This period encompasses significant shifts in public sentiment driven by
evolving pandemic conditions, policy changes, and vaccine rollout. Each
tweet was annotated with a compound sentiment score ranging from $-1$ (most
negative) to $+1$ (most positive) using the Valence Aware Dictionary and
sEntiment Reasoner (VADER) \cite{HuttoGilbert2014}, a lexicon and rule-based
sentiment analysis model optimized for social media text.


The main contributions of this paper are as follows. We provide the first
application of convexification-based MFG forecasting to real-world dynamic sentiment
data spanning an extended time period. Furthermore, we demonstrate that the
solution obtained via convexification not only matches observed public sentiment data but
also satisfies the governing MFG equations with small residuals when
appropriately parameterized. To our knowledge, this work provides the first
empirical evidence that MFGs can effectively govern complex temporal
patterns in public sentiment, including sudden shifts and gradual trends.

While the current study relies on manual parameter calibration due to
the absence of systematic coefficient identification methods, our results
provide a strong proof-of-concept evidence for MFG-based forecasting. We
emphasize that parameter calibration is necessary here because the dataset
provides no information about the true model coefficients, and developing
systematic methods to identify them from observed sentiment data remains an
important open problem. Future work will focus on developing robust
solvers for coefficient inverse problems to determine these coefficients
directly from the measured data--an objective that appears feasible given
the successful application of convexification to several coefficient inverse
problems in MFG contexts \cite{KL,MFGbook,MFGCAMWA,MFG8IPI}.

The remainder of this paper is organized as follows. Section \ref%
{sec:problem_statement} formulates the MFG system and the forecasting
problem. Section \ref{sec:theoretical_justification} summarizes the
theoretical justification for the convexification method. Section \ref%
{sec:dataset_description} provides details and insights regarding the
sentiment dataset. Section \ref{sec:forecasting_methodology} describes the
forecasting procedure and parameter calibration strategy. Section \ref%
{sec:numerical_results} presents numerical results for ten consecutive time
periods. Section \ref{sec:discussion_future_directions} discusses our
findings and potential future directions.

\section{The MFG System and Problem Statement}

\label{sec:problem_statement}

Below $\Omega \subset \mathbb{R}^{n}$ is a bounded domain with the piecewise
smooth boundary $\partial \Omega .$ We consider the following $n$%
-dimensional MFG system \cite{Achdou}: 
\begin{gather}
u_{t}+\beta \Delta u+\frac{1}{2}r|\nabla u|^{2}+\int_{\Omega
}K(x,y)m(y,t)dy=0,\quad (x,t)\in Q_{T},  \label{2.1} \\
m_{t}-\beta \Delta m+\text{div}\left( rm\nabla u\right) =0,\quad (x,t)\in
Q_{T},  \label{2.2}
\end{gather}%
where $T>0$ is a constant, $\beta =\beta (x,t)$ is the diffusion
coefficient, $r=r(x,t)$ is the drift coefficient, $K(x,y)$ is the
interaction kernel, and 
\begin{equation}
Q_{T}=\Omega \times (0,T).  \label{2.20}
\end{equation}%
Here $\Omega $ is viewed as a state space, and the system may be interpreted
as follows. We consider a population of continuum-many rational agents. For $%
x\in \Omega $ and $t\in \lbrack 0,T]$, we view $m(x,t)$ as the density of
the population occupying state $x$ at the time $t$, and we view $u(x,t)$ as
the value (or expected payoff) of an agent being at state $x$ at time $t$.
In the MFG system, \eqref{2.1} is the Hamilton-Jacobi-Bellman (HJB)
equation, and \eqref{2.2} is the Fokker-Planck-Kolmogorov (FPK) equation.
{The conventional MFG system \cite{LasryLions2007} is equipped with an initial condition for $m$ and a terminal condition for $u$: 
\begin{equation*}
    u(x,T) = u_T(x), \quad m(x,0) = m_0(x), \quad x \in \Omega.
\end{equation*}
The HJB equation is viewed backwards in time while the FPK equation is viewed forwards in time. We note that uniqueness for this system has only been established under restrictive monotonicity assumptions \cite{LasryLions2007}.
}

{
In this paper we are interested in a forecasting problem in which we only have access to initial data for $m$ and we would nonetheless like to make projections about $m(x,t)$ for $t > 0$. Thus, although we are modeling the function $m(x,t)$ as part of a solution to a conventional MFG system, we do not have access to the value function $u(x,t)$ at any time, including the terminal time. It is well-known that this problem is unstable, i.e., ill-posed. Nonetheless, recent work has established that if the initial condition of the value function $u(x,0) = u_0(x)$ is prescribed, then the problem has H\"{o}lder stability and uniqueness \cite{MFGbook,MFGSAM}. (Note that similar results are not available if we assume instead that the terminal condition $u(x,T) = u_T(x)$ is prescribed.) Thus, in this work we treat the initial condition $u(x,0)$ of the value function $u(x,t)$ as a latent parameter in the problem and develop our numerical method assuming that this parameter is known. This naturally leads us to consider the following initial and boundary conditions:
\begin{gather}
u(x,0)=u_{0}(x),\quad m(x,0)=m_{0}(x),\quad x\in \Omega ,  \label{2.3} \\
\frac{\partial u}{\partial \nu }=\frac{\partial m}{\partial \nu }=0\quad 
\text{on }\partial \Omega ,\quad \forall t\in (0,T),  \label{2.4}
\end{gather}%
where $\nu $ is the outward unit normal vector on the boundary $\partial
\Omega $.
}

\begin{remark}
We refer to our previous work \cite{Kliba2025} for a detailed discussion on the challenges of having both initial conditions for $u$ and $m$. For example, in the case where $\beta$ is a positive constant, the term $u_t + \beta \Delta u$ combined with initial condition $u(x,0) = u_0(x)$ introduces extreme instability into the task of solving for $u$, rendering standard time-marching numerical approaches ineffective. {The convexification method proposed here offers a balance between accuracy of the solution and stability of the problem up to a certain point in time.} 
\end{remark}

{
\begin{remark}
 There is a trade-off associated with treating the initial value of the value function as a latent parameter. Indeed, this formulation means that we must try to use the available data to identify this parameter. See Section \ref{sec:initial_value_estimation} for our approach to this aspect of the problem.
\end{remark}
}



Now consider the Sobolev space 
\begin{equation*}
H_0^{2}(Q_T) := \left\{ f \in H^{2}(Q_T): \frac{\partial f}{\partial \nu} =
0 \quad \text{on } \partial \Omega,\ \forall t \in (0,T) \right\}, 
\end{equation*}
where $H^2(Q_T)$ denotes the standard Sobolov space of functions on $Q_T$
with square-integrable weak derivatives up to order two. We formulate the
forecasting problem for the MFG system \eqref{2.1}--\eqref{2.4} as follows.

\paragraph{The Forecasting Problem}

Given the functions $r$, $\beta $, $K$, and the initial conditions $m_{0}(x)$
and $u_{0}(x)$, find $u,m\in H_{0}^{2}(Q_{T})$ that satisfy the MFG system %
\eqref{2.1}--\eqref{2.4} for all times $t\in (0,t_{0})$, where $t_{0}$ is a
point of time between $0$ and $T$.

H\"{o}lder stability and uniqueness of the solution to the forecasting
problem have been established in \cite{MFGbook} and \cite{MFGSAM}. However, solving the MFG system \eqref{2.1}--%
\eqref{2.4} numerically presents significant challenges. First, the HJB
equation with initial condition \eqref{2.3} exhibits a similar behavior as
the backward heat equation, meaning that solutions typically grow
unboundedly as time advances, and conventional forward time-marching methods
are therefore inapplicable. Second, the inherent nonlinearity of the coupled
system renders standard optimization approaches non-convex, making them
susceptible to local minima. 

To overcome these difficulties, we employ the convexification method
developed in our previous work \cite{Kliba2025}. This approach transforms
the original problem into a strongly convex optimization problem through the
introduction of an appropriate CWF. The resulting convexified functional
admits a unique global minimizer, which well approximates the solution of
the MFG system while maintaining stability. The theoretical justification
and rigorous convergence analysis of the convexification method are
presented in detail in \cite{Kliba2025}. For completeness and reader
convenience, we summarize the key theoretical results for the
one-dimensional case in the next section.

\section{Theoretical Justification of Convexification for the MFG System}

\label{sec:theoretical_justification}

We focus on the case where $\Omega = (-1,1) \subset \mathbb{R}$. This
setting is particularly well-suited for our application to public sentiment
forecasting, as sentiment scores naturally lie on a one-dimensional spectrum
ranging from negative to positive values. The MFG system \eqref{2.1}--%
\eqref{2.4} reduces to 
\begin{gather}
u_{t} + \beta u_{xx} + \frac{1}{2}ru_x^2 + \int_{\Omega} K(x,y) m(y,t) dy =
0, \quad (x,t) \in Q_T,  \label{3.01} \\
m_{t}-\beta m_{xx}+\partial _{x}\left( r mu_{x}\right) = 0, \quad (x,t) \in
Q_T,  \label{3.02} \\
u(x,0) = u_0(x), \quad m(x,0) = m_0(x), \quad x \in (-1,1),  \label{3.03} \\
u_x(-1,t) = u_x(1,t) = m_x(-1,t) = m_x(1,t) = 0, \quad t \in (0,T).
\label{3.04}
\end{gather}

For the sake of completeness, we summarize below the key theoretical results
from \cite{Kliba2025} that justify the convexification method for solving
the forecasting problem for \eqref{3.01}--\eqref{3.04}. We begin with the
following assumption regarding the coefficients $\beta $ and $r$, as well as
the kernel $K$.

\begin{assumption}
\label{assumption:1} Assume that 
\begin{gather*}
\beta =const.>0,r\in C^{1}(\overline{Q_{T}}),\Vert r\Vert _{C^{1}(\overline{%
Q_{T}})}\leq M, \\
K\in L^{\infty }(\Omega \times \Omega ),\quad \Vert K\Vert _{L^{\infty
}(\Omega \times \Omega )}\leq M,
\end{gather*}%
for some constant $M>0$.
\end{assumption}

Let $c>2$ be such that 
\begin{equation}
\frac{c^{2}}{T+c}\geq 2,  \label{eq:c_condition}
\end{equation}%
and define the CWF as 
\begin{equation*}
\phi _{\lambda }(t):=e^{(T-t+c)^{\lambda }},\quad t\in (0,T),
\end{equation*}%
for some $\lambda >1$. The following two Carleman estimates play a crucial
role in the convexification method.

\begin{theorem}
\label{theorem:carleman_1} (\cite[Section 2.3.1]{MFGbook})
Let $c$ be the constant satisfying \eqref{eq:c_condition}. Then, there
exists a sufficiently large number $\lambda _{0,1}=\lambda _{0,1}(\beta
,c,M,T)\geq 1$ such that for all $\lambda \geq \lambda _{0,1}$ and for all
functions $u\in H_{0}^{2}(Q_{T})$, the following Carleman estimate holds: 
\begin{multline*}
\int_{Q_{T}}(u_{t}+\beta u_{xx})^{2}\phi _{\lambda }^{2}\ dxdt\geq  \\
\geq C_{1}\sqrt{\lambda }\int_{Q_{T}}u_{x}^{2}\phi _{\lambda }^{2}\
dxdt+C_{1}\lambda ^{2}\int_{Q_{T}}u^{2}\phi _{\lambda }^{2}\ dxdt- \\
-C_{1}e^{2c\lambda }\int_{\Omega }(u_{x}^{2}+u^{2})(x,T)\ dx-C_{1}\lambda
(T+c)^{\lambda }e^{2(T+c)^{\lambda }}\int_{\Omega }u^{2}(x,0)\ dx,
\end{multline*}%
where the constant $C_{1}=C_{1}(\beta ,c,M,T)>0$ depends only on the listed
parameters.
\end{theorem}

\begin{theorem}
\label{theorem:carleman_2} (\cite[Section 2.3.2]{MFGbook})
Let $c$ be the constant satisfying \eqref{eq:c_condition}. Then, there
exists a sufficiently large number $\lambda _{0,2}=\lambda _{0,2}(\beta
,c,M,T)\geq 1$ such that for all $\lambda \geq \lambda _{0,2}$ and for all
functions $u,v\in H_{0}^{2}(Q_{T})$, the following quasi-Carleman estimate
holds: 
\begin{multline*}
\int_{Q_{T}}(u_{t}-\beta u_{xx}+rv_{xx})^{2}\phi _{\lambda }^{2}\ dxdt\geq 
\\
\geq \lambda c^{\lambda -1}\int_{Q_{T}}u_{x}^{2}\phi _{\lambda }^{2}\ dxdt+%
\frac{\lambda ^{2}}{4}c^{2\lambda -2}\int_{Q_{T}}u^{2}\phi _{\lambda }^{2}\
dxdt- \\
-C_{2}\lambda (T+c)^{\lambda }\int_{Q_{T}}v_{x}^{2}\phi _{\lambda }^{2}\
dxdt-C_{2}\lambda (T+c)^{\lambda }e^{2(T+c)^{\lambda }}\int_{\Omega
}u^{2}(x,0)\ dx,
\end{multline*}%
where the constant $C_{2}=C_{2}(\beta ,c,M,T)>0$ depends only on the listed
parameters.
\end{theorem}

In \cite[Section 2.3.2]{Kliba2025}, we need our functions $u,m\in
C^{3}\left( \overline{Q_{T}}\right) .$ Hence, by Sobolev embedding theorem,
we use $s\geq 5$ in the following function spaces: 
\begin{gather*}
H_{2}^{s}(Q_{T}):=\left\{ 
\begin{array}{c}
(u,m)\in H^{s}(Q_{T})\times H^{s}(Q_{T}): \\ 
\Vert (u,m)\Vert _{H_{2}^{s}(Q_{T})}^{2}:=\Vert u\Vert
_{H^{s}(Q_{T})}^{2}+\Vert m\Vert _{H^{s}(Q_{T})}^{2}<\infty 
\end{array}%
\right\} , \\
H_{2,0}^{s}(Q_{T}):=\left\{ 
\begin{array}{ll}
(u,m)\in H_{2}^{s}(Q_{T}): & u_{x}(-1,\cdot )=u_{x}(1,\cdot )=0, \\ 
& m_{x}(-1,\cdot )=m_{x}(1,\cdot )=0,%
\end{array}%
\right\} , \\
H_{2,0}^{s}(\Omega ):=\left\{ 
\begin{array}{c}
(f,g): \\ 
\Vert (f,g)\Vert _{H_{2}^{s}(\Omega )}^{2}:=\Vert f\Vert _{H^{s}(\Omega
)}^{2}+\Vert g\Vert _{H^{s}(\Omega )}^{2}<\infty , \\ 
f_{x}(-1)=f_{x}(1)=g_{x}(-1)=g_{x}(1)=0%
\end{array}%
\right\} , \\
H_{2,0,0}^{s}(Q_{T}):=\left\{ (u,m)\in H_{2,0}^{s}(Q_{T}):u(\cdot
,0)=m(\cdot ,0)=0\right\} ,
\end{gather*}%
and denote the scalar product in $H_{2}^{s}(Q_{T})$ by $[\cdot ,\cdot ]$. In
addition, for 
\begin{equation}
\gamma \in (0,1),  \label{eq:gamma_condition}
\end{equation}%
we let $Q_{\gamma T}:=\Omega \times (0,\gamma T)$ and define the space 
\begin{equation*}
H^{1,0}(Q_{\gamma T}):=\left\{ u:\Vert u\Vert _{H^{1,0}(Q_{\gamma
T})}^{2}:=\Vert u_{x}\Vert _{L^{2}(Q_{\gamma T})}^{2}+\Vert u\Vert
_{L^{2}(Q_{\gamma T})}^{2}<\infty \right\} .
\end{equation*}%
Let $R>0$ be an arbitrary number. We assume that the initial conditions
satisfy 
\begin{equation}
(u_{0},m_{0})\in H_{2,0}^{s}(\Omega ),\quad \Vert (u_{0},m_{0})\Vert
_{H_{2}^{s}(\Omega )}<R,  \label{eq:initial_conditions}
\end{equation}%
and define the following set of admissible solutions 
\begin{equation*}
B(R):=\left\{ 
\begin{array}{c}
(u,m)\in H_{2,0}^{s}(Q_{T}):\Vert (u,m)\Vert _{H_{2}^{s}(Q_{T})}<R, \\ 
u(x,0)=u_{0}(x),\quad m(x,0)=m_{0}(x),\quad x\in \Omega 
\end{array}%
\right\} .
\end{equation*}%
Let $L_{1}(u,m)$ and $L_{2}(u,m)$ be the operators on the left-hand side of %
\eqref{3.01} and \eqref{3.02}, respectively, and define 
\begin{equation*}
q=q_{\lambda }(c,T):=\frac{1}{\lambda (T+c)^{\lambda -1}}.
\end{equation*}%
The convexification numerical algorithm aims to find the solution of the MFG
system (\ref{3.01})-(\ref{3.02}) with initial and boundary conditions (\ref%
{3.03}), (\ref{3.04}) by minimizing the weighted functional $J_{\lambda
,\alpha }(m,u):\overline{B(R)}\rightarrow \mathbb{R}$ defined by 
\begin{multline}
J_{\lambda ,\alpha }\left( m,u\right) :=e^{-2ac^{\lambda
}}\int\limits_{Q_{T}}\left[ L_{1}(u,m)^{2}+qdL_{2}(u,m)^{2}\right] \phi
_{\lambda }^{2}\ dxdt+  \label{3.1} \\
+\alpha \Vert (u,m)\Vert _{H_{2}^{s}(Q_{T})}^{2},
\end{multline}%
where $\alpha \in (0,1)$ is the regularization parameter and $a,d>0$ are two
numbers to be chosen for the numerical implementation. The constant $%
e^{-2ac^{\lambda }}$ is introduced to partially balance the two terms in %
\eqref{3.1} since the maximum value of the CWF is 
\begin{equation*}
\max_{t\in \lbrack 0,T]}\phi _{\lambda }(t)=\phi _{\lambda
}(0)=e^{(T+c)^{\lambda }}.
\end{equation*}%
The Carleman estimates of Theorems \ref{theorem:carleman_1} and \ref%
{theorem:carleman_2} allow us to establish the following convexity and error
estimates, which are the main theoretical results of \cite{Kliba2025}.

\begin{theorem}
If Assumption \ref{assumption:1} along with \eqref{eq:c_condition} and %
\eqref{eq:initial_conditions} hold, then: 

\begin{enumerate}
\item The Fr\'{e}chet derivative $J_{\lambda ,\alpha }^{\prime }\in
H_{2,0,0}^{s}(Q_{T})$ of the functional $J_{\lambda ,\alpha }$ exists at
each point $(u,m)\in \overline{B(R)}$ and is Lipschitz continuous on $%
\overline{B(R)}$. 

\item There exists a sufficiently large number $\lambda _{1}=\lambda
_{1}(\beta ,M,R,c,a,d,T)\geq \max \{\lambda _{0,1},\lambda _{0,2}\}$ such
that for any $\lambda \geq \lambda _{1}$ and for any $\alpha \in \lbrack
2e^{-\lambda c^{\lambda }},1)$, the functional $J_{\lambda ,\alpha }$ is
strongly convex on the set $\overline{B(R)}$, i.e. there exists a constant $%
C=C(\beta ,M,R,c,a,d,\gamma ,T)>0$ such that 
\begin{multline}
J_{\lambda ,\alpha }(u_{2},m_{2})-J_{\lambda ,\alpha
}(u_{1},m_{1})-[J_{\lambda ,\alpha }^{\prime
}(u_{1},m_{1}),(u_{2}-u_{1},m_{2}-m_{1})]\geq  \\
\geq C\left( \Vert u_{2}-u_{1}\Vert _{H_{1,0}^{s}(Q_{\gamma T})}^{2}+\Vert
m_{2}-m_{1}\Vert _{H_{1,0}^{s}(Q_{\gamma T})}^{2}\right) + \\
+\frac{\alpha }{2}\Vert (u_{2}-u_{1},m_{2}-m_{1})\Vert
_{H_{2}^{s}(Q_{T})}^{2}, \\
\forall (u_{1},m_{1}),(u_{2},m_{2})\in \overline{B(R)},\quad \lambda \geq
\lambda _{1}.
\end{multline}%
Numbers $\lambda _{1}$ and $C$ depend only on listed parameters.

\item If $\lambda \geq \lambda _{1}$ and $\alpha \in \lbrack 2e^{-\lambda
c^{\lambda }},1)$, then there exists a unique minimizer $(u_{\min ,\lambda
,\alpha },m_{\min ,\lambda ,\alpha })\in \overline{B(R)}$ of the functional $%
J_{\lambda ,\alpha }$ on the set $\overline{B(R)}$. Moreover, the following
estimate holds for all $(u,m)\in \overline{B(R)}$: 
\begin{equation}
\left[ J_{\lambda ,\alpha }^{\prime }(u_{\min ,\lambda ,\alpha },m_{\min
,\lambda ,\alpha }),(u_{\min ,\lambda ,\alpha }-u,m_{\min ,\lambda ,\alpha
}-m)\right] \leq 0.
\end{equation}
\end{enumerate}
\end{theorem}

It is always assumed in the theory of
ill--posed problems that there exists an ideal or true solution of such a
problem with the ideal noiseless data \cite{Tikhonov1995}. Hence, we formulate now the
following accuracy estimate of \cite[Section 2.3.2]{Kliba2025}:

\begin{theorem}
Assume that \eqref{eq:gamma_condition}, \eqref{eq:initial_conditions} hold
and that there exists an ideal solution $(u^{\ast },m^{\ast })$ of the MFG
system \eqref{2.1}--\eqref{2.4} with the exact initial conditions $%
(u_{0}^{\ast },m_{0}^{\ast })\in H_{2,0}^{s}(\Omega )$ such that 
\begin{equation*}
(u^{\ast },m^{\ast })\in B^{\ast }(R):=\left\{ 
\begin{array}{c}
(u,m)\in H_{2,0}^{s}(Q_{T}):\Vert (u,m)\Vert _{H_{2}^{s}(Q_{T})}<R, \\ 
u(x,0)=u_{0}^{\ast }(x),\ m(x,0)=m_{0}^{\ast }(x)%
\end{array}%
\right\} .
\end{equation*}%
Let $\delta \in (0,1)$ be the level of noise in the initial value $%
(u_{0},m_{0})$, i.e. 
\begin{equation}
\Vert u_{0}-u_{0}^{\ast }\Vert _{H^{2}(\Omega )}<\delta ,\quad \Vert
m_{0}-m_{0}^{\ast }\Vert _{H^{2}(\Omega )}<\delta .
\end{equation}%
Then, there exists $\delta _{0}=\delta _{0}(M,R,c,a,d,\gamma ,T)\in (0,1)$
such that for any $\delta \in (0,\delta _{0})$ and for 
\begin{gather*}
\lambda =\lambda (\delta ):=\frac{1}{2\ln (T+c)}\ln \left[ \left( \ln \left(
\delta ^{-1/3}\right) \right) \right] , \\
\alpha =\alpha (\delta ):=2e^{-(a-1)c^{\lambda (\delta )}},
\end{gather*}%
the unique minimizer $(u_{\min },m_{\min })\in \overline{B(R)}$ of the
functional $J_{\lambda (\delta ),\alpha (\delta )}$ with the noisy initial
data $(u_{0},m_{0})$ satisfies the following estimate: 
\begin{equation*}
\Vert u_{\min }-u^{\ast }\Vert _{H^{1,0}(Q_{\gamma T})}+\Vert m_{\min
}-m^{\ast }\Vert _{H^{1,0}(Q_{\gamma T})}\leq C_{1}\sqrt{\delta },
\end{equation*}%
where the constant $C_{1}=C_{1}(M,R,c,a,d,\gamma ,T)>0$ only depends on the
listed parameters.
\end{theorem}

In summary, the above theorems guarantee the existence and uniqueness of the
minimizer of the functional $J_{\lambda ,\alpha }$ as well as the
convergence of this minimizer to the true solution of the MFG system as the
noise in the initial data tends to zero. These theoretical results provide
the mathematical foundation for our forecasting approach.

We now turn from the theoretical framework to its practical application to
real-world sentiment data drawn from public COVID-19 discussions. In the
next section, we describe the characteristics of the dataset used in our
study.

\section{Dataset Description}

\label{sec:dataset_description}

Our study is based on a real dataset of highly engaged tweets related to
general public discourses on COVID-19, collected using the Brandwatch
platform\footnote{\url{https://www.brandwatch.com/}}, with a sampling rate
of 1\%. We focused on tweets with an engagement score of 10 or higher,
ensuring that the dataset represents content that elicited notable public
attention. The engagement score on X is calculated as the sum of a tweet's
likes, reposts (retweets), and replies. On Brandwatch, the engagement score
is a platform-specific aggregate, i.e., each social network platform
combines its own relevant interaction indicators into a single ``sum-all'' value. To gather relevant content, we used
a broad keyword-based query encompassing medical, colloquial, and
politicized references to COVID-19 as is shown in Table \ref%
{tab:covid_keywords}. The keywords were searched in either the bodies or the
titles of the tweets, and results were filtered to include only those from X.

\begin{table}[htbp]
\centering
\resizebox{\textwidth}{!}{\begin{tabular}{@{}>{\raggedright\arraybackslash}p{0.25\textwidth}                >{\raggedright\arraybackslash}p{0.7\textwidth}@{}}
\toprule
\multicolumn{1}{c}{\textbf{Category}} & 
\multicolumn{1}{c}{\textbf{Keywords}} \\ 
\midrule
\small \textit{General Keywords} & \texttt{ncov, ncov-19, sars, SARS-CoV-2, coronavirus, pandemic, pheic} \\

\small \textit{Politicized Terms} & \texttt{"wuhan virus", "china virus", "wuhan pneumonia", "wuhan flu", kungflue} \\

\small \textit{Standard Terms of COVID-19} & \texttt{covid19, "covid-19", covid, "covid 19"} \\

\small \textit{Official Terminology} & \texttt{"Public Health Emergency of International Concern"} \\
\bottomrule
\end{tabular}}
\caption{COVID-19 Keywords Used in Query}
\label{tab:covid_keywords}
\end{table}

The timeframe of this real-world dataset spans from March 2, 2020, at 00:01
a.m. to April 10, 2022, at 11:59 p.m. (UTC), comprising a total of 47,181
mentions. This period encompasses several pivotal phases of the pandemic,
including the initial outbreak response, the implementation of lockdown
measures, and the subsequent introduction of travel bans, school closures,
and vaccination roll-out campaigns, among other major public health
interventions. The selected timeframe thus captures the dynamic evolution of
public sentiment as the pandemic unfolded.


Each tweeted mention of the public discussion on COVID-19 was assigned a
sentiment score (compound score) ranging from $-1$ (most negative) to $+1$
(most positive) using the VADER \cite{HuttoGilbert2014} sentiment analysis
tool. VADER is a rule-based sentiment analysis tool specifically designed
for social media text, as it accounts for factors such as sentiment
intensity, capitalization, punctuation, and the presence of emoticons.
Weekly aggregations of these scores were subsequently performed to derive
temporal sentiment distributions, enabling the analysis of shifts in public
sentiment and emotional expressions over time.


\begin{figure}[ht!]
\centering
\includegraphics[width=1\linewidth]{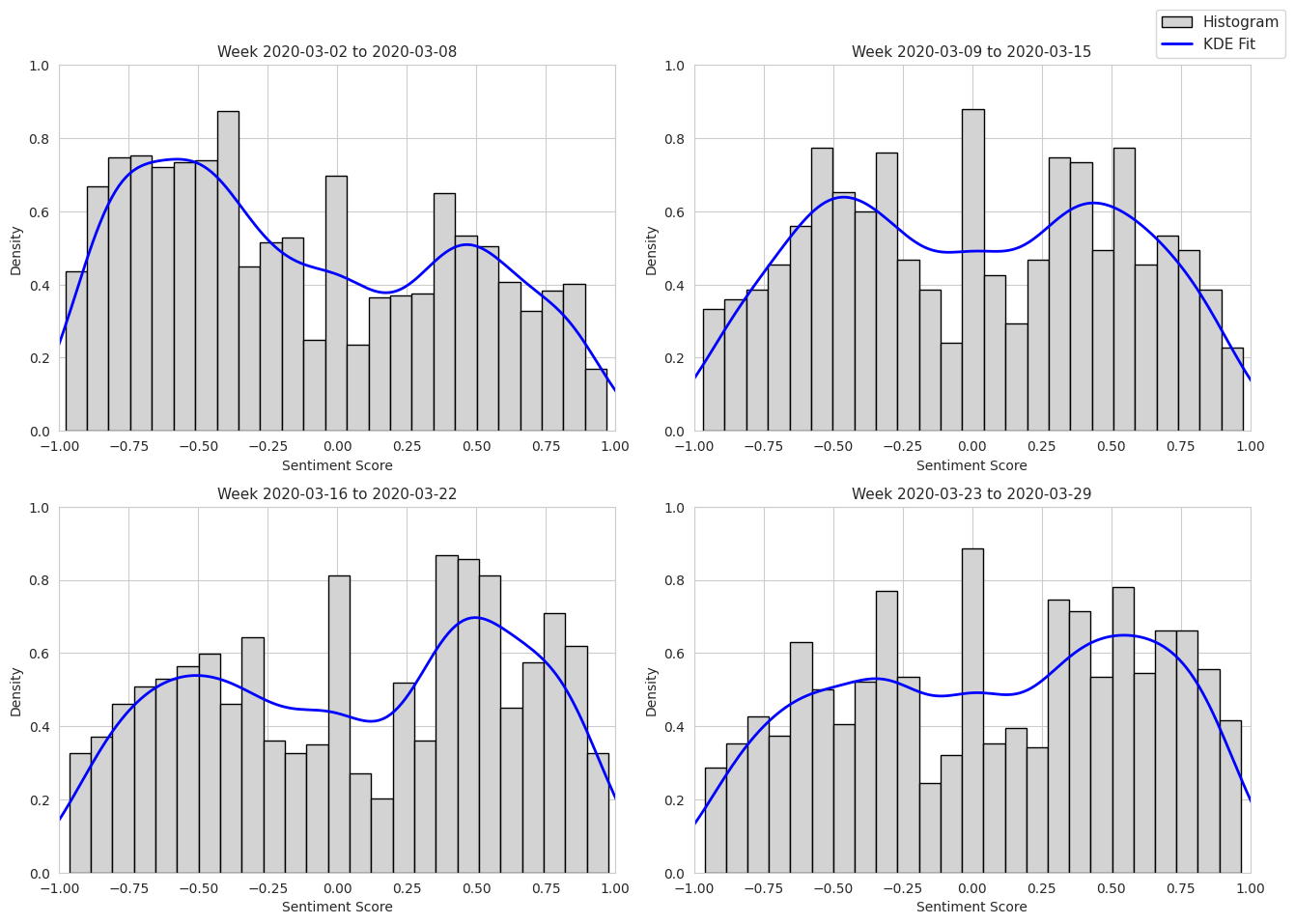}  
\caption{Sentiment distributions of public discussions on COVID-19 across
four consecutive weeks from March~2 to March~29,~2020. Each subplot displays
the histogram of VADER compound sentiment scores for the corresponding week,
overlaid with a KDE fit curve shown in solid blue. The distributions reveal
temporal variations in sentiment polarity and dispersion during the early
stage of the pandemic. Note: the histograms have been normalized to represent the density of sentiment scores in each of the four weeks, not the original counts.}
\label{fig:sentiment_density_first_four_weeks}
\end{figure}


The sentiment probability density $m(x,t)$ was estimated from discrete
sentiment scores using kernel density estimation (KDE) with a Gaussian
kernel. Figure \ref{fig:sentiment_density_first_four_weeks} presents
representative histograms of sentiment scores and their corresponding
estimated densities for the first four weeks of the dataset. A salient
feature of the data is the strong concentration of sentiment scores around
the neutral point ($0$), producing a pronounced peak in the density across
all subplots in Figure \ref{fig:sentiment_density_first_four_weeks}. This
observation likely reflects the communicative characteristics of social
media interactions, where a large proportion of posts contain neutral or
context-independent content, or express mixed and ambivalent attitudes
toward the unfolding pandemic situation.

\section{Forecasting Methodology}

\label{sec:forecasting_methodology}

\subsection{Parameters of the Minimization Functional}

We fix the parameters $K(x,y)$, $\lambda $, $a$, $c$, $d$, and $\alpha $ in
the minimization functional \eqref{3.1}. Given the long two-year span of the
dataset and substantial variations in public sentiment driven by evolving
pandemic conditions, it is unlikely that other parameters would
remain the same for the entire two years period. Hence, we partition the
data into smaller temporal intervals for separate analysis. The dataset
contains 110 weeks of sentiment densities from March 2, 2020 to April 10,
2022. We divide this into 10 periods of 11 weeks each, and treat each period
as an independent forecasting problem. For each period, the first week is
associated with $t=0$, and the eleventh week with $t=T=1$.

According to the problem formulation in Section \ref{sec:problem_statement},
forecasting requires knowledge of the coefficients $\beta $, $r(x,t)$, and
the initial conditions $m(x,0)$, $u(x,0)$. 
As described in Section \ref{sec:dataset_description}, the initial density
$m(x,0)$ is obtained from the weekly sentiment histogram using Gaussian KDE. In fact, the dataset produces an \textit{observed} density $m(x,t)$ at every weekly time point used in our experiments.
However, the diffusion
coefficient $\beta $, drift coefficient $r$, and initial value function $%
u(x,0)$ are unknown. This lack of coefficient information presents a
significant challenge. In the absence of systematic identification methods,
we conducted extensive numerical experiments to experimentally calibrate
these missing values. First, we describe our approach to estimating the
initial value function $u(x,0)$ from the data, and then we describe our
experimental strategy to calibrate $\beta $ and $r$.

\subsection{Initial Value Function Estimation}

\label{sec:initial_value_estimation}

We start with the FPK equation and the Neumann boundary conditions: 
\begin{gather}
m_{t}-\beta m_{xx}+r\partial _{x}\left( mu_{x}\right) =0,x\in \left(
-1,1\right) ,t\in \left( 0,T\right) \label{1} \\
u_{x}\left( -1,t\right) =u_{x}\left( 1,t\right) =0, \\
m_{x}\left( -1,t\right) =m_{x}\left( 1,t\right) =0.
\end{gather}%
Let $h$ be the temporal step size. We can approximate $m_{t}\left(
x,0\right) $ using the following second-order finite difference scheme: 
\begin{equation}
m_{t}\left( x,0\right) \approx \frac{-3m\left( x,0\right) +4m(x,h)-2m\left(
x,2h\right) }{2h}.  \label{5.1}
\end{equation}
Here we have used the sentiment densities estimated from the first three
weeks of data, i.e., $m(x,0)$, $m(x,h)$, and $m(x,2h)$. This is the
trade-off we make to achieve a more accurate approximation of $m_{t}(x,0)$
while minimizing the number of weeks used for initial value estimation.

Assuming that $r(x,0)m(x,0)\neq 0$ for all $x\in \left[ -1,1\right] $, by (%
\ref{1}), we have 
\begin{equation}
u_{xx}\left( x,0\right) +\frac{m_{x}\left( x,0\right) }{m\left( x,0\right) }%
u_{x}\left( x,0\right) =p\left( x,0\right) ,  \label{8}
\end{equation}%
where 
\begin{equation*}
p\left( x,0\right) :=\frac{1}{rm\left( x,0\right) }\left[ \beta m_{xx}\left(
x,0\right) -m_{t}\left( x,0\right) \right] .
\end{equation*}%
Denote 
\begin{equation}
u_{x}\left( x,0\right) =v\left( x,0\right) ,  \label{10}
\end{equation}
with $v(-1,0)=0$, so that (\ref{8}) becomes%
\begin{equation}
\left. v_{x}\left( x,0\right) +\frac{m_{x}\left( x,0\right) }{m\left(
x,0\right) }v\left( x,0\right) =p\left( x,0\right) ,\right.   \label{11}
\end{equation}%
\begin{equation}
v(-1,0)=0.  \label{110}
\end{equation}%
Solving the initial value problem (\ref{10})-(\ref{110}) via an explicit
formula gives 
\begin{equation}
u_{x}\left( x,0\right) =\int\limits_{-1}^{x}p\left( y,0\right) \exp \left(
\int\limits_{y}^{x}\frac{m_{x}\left( s,0\right) }{m\left( s,0\right) }%
ds\right) .  \label{13}
\end{equation}

{It follows from (\ref{10})-(\ref{110}) that it not necessary that
the second boundary condition }%
\begin{equation}
u_{x}\left( 1,0\right) =0  \label{14}
\end{equation}%
{would be satisfied for the function }$u_{x}\left( x,0\right) $%
{\ given by (\ref{13}). Thus, since we deal with approximations, we
need to ensure somehow that   ensure the condition (\ref{13}) is still
satisfied. Hence, to ensure the latter, }we consider a smooth cut-off
function $\chi \left( x\right) \in C^{2}\left[ -1,1\right] $ such that 
\begin{equation*}
\chi \left( x\right) =\left\{ 
\begin{array}{ll}
1, & x\in (-1,1-\varepsilon ] \\ 
1-3\left( \frac{x-1+\varepsilon }{\varepsilon }\right) ^{2}+2\left( \frac{%
x-1+\varepsilon }{\varepsilon }\right) ^{3}, & x\in (1-\varepsilon ,1)%
\end{array}%
\right. 
\end{equation*}%
for a small $\varepsilon \in \left( 0,1\right) $. Then replacing (\ref{13})
with 
\begin{equation*}
u_{x}\left( x,0\right) \approx \chi \left( x\right)
\int\limits_{-1}^{x}p\left( y,0\right) \exp \left( \int\limits_{y}^{x}\frac{%
m_{x}\left( s,0\right) }{m\left( s,0\right) }ds\right) ,
\end{equation*}%
we obtain 
\begin{equation*}
u\left( x,0\right) \approx u\left( -1,0\right) +\int\limits_{-1}^{x}\chi
\left( z\right) \left[ \int\limits_{-1}^{z}p\left( y,0\right) \exp \left(
\int\limits_{y}^{s}\frac{m_{x}\left( s,0\right) }{m\left( s,0\right) }%
ds\right) dy\right] dz.
\end{equation*}

\begin{remark}
The above process also illustrates the necessity of additional steps to
ensure that $u(x,0)$ satisfies the Neumann boundary condition at both
endpoints $x=-1$ and $x=1$. Specifically, we must employ a cut-off function $%
\chi (x)$. This highlights the extreme ill-posedness of the forecasting
problem, indicating that the existence of the solution cannot be guaranteed
unless some restrictive and yet unknown conditions are imposed.
\end{remark}

Since the available data do not provide information about the boundary value 
$u(-1,0)$, we treat it as a free parameter to be calibrated alongside the
coefficients $\beta $ and $r$ (see Section \ref{sec:calibration_procedure}).
Note, again, that for each 11-weeks period, the estimation of the initial
condition $u(x,0)$ requires sentiment densities from the first three weeks,
namely $m(x,0)$, $m(x,h)$, and $m(x,2h)$, due to \eqref{5.1}. Consequently,
the effective forecasting horizon for each period spans weeks 4--11,
comprising eight weeks.

Next, we describe our experimental procedure to calibrate the unknown
parameters $\beta$, $r$, and $u(-1,0)$.

\subsection{Calibration Procedure}

\label{sec:calibration_procedure}

For each fixed 11-week period, we assume that $\beta $ and $r$ remain
constant. Recall that, as mentioned above, the observed density $m(x,t)$ is available at every weekly time point $t$. We perform the following experiment:

\begin{enumerate}
\item[S1.] Initialize parameters $\beta $, $r$, and $u(-1,0)$ with trial
values. The trial values are chosen randomly for the first period. For
subsequent periods, we use the calibrated values from the previous period as initial guesses.

\item[S2.] Estimate the initial value function $u(x,0)$ following the
procedure described in Section \ref{sec:initial_value_estimation}. 

\item[S3.] Using the observed $m(x,0),m(x,h)$, and $m(x,2h)$ from the first three weeks, the estimated $u(x,0)$, and the current values of $\beta $ and $r$, solve the forecasting
problem via convexification to predict sentiment densities weeks 4--11, comprising eight weeks. 

\item[S4.] Evaluate the match between the convexification solution and
observed sentiment densities using visual inspection. This step
essentially verifies how well we chose parameters in S3.

\item[S5.] Adjust $u(-1,0)$, $\beta$, and $r$ and repeat S3--S5 until
satisfactory agreement is achieved.
\end{enumerate}

\begin{remark}
Although this procedure does not constitute a true predictive test in the
strict sense (since we observe the outcomes during our calibration
procedure), it serves as a proof-of-concept demonstration that the MFG model
can be parameterized to reproduce observed sentiment dynamics. From this
procedure, we are able to draw conclusions of the following type. Given an
appropriate choice of model parameters $\beta $ and $r$ and initial
measurements $m(x,0)$ and $u(x,0)$, public sentiments in the dataset can be
governed well by the proposed MFG model. Furthermore, the convexification
method provides a systematic and stable approach to compute these forecasts.
We note that our future work will focus on developing inverse problem
techniques to determine $\beta $ and $r$ directly from data, enabling
genuine predictive capability. See Section \ref%
{sec:discussion_future_directions} for further discussion.
\end{remark}

\section{Numerical Results}

\label{sec:numerical_results}

\subsection{Implementation Details}

\label{sec:Implementation}

Forecasting is performed by minimizing a slightly modified version of the
functional \eqref{3.1}: 
\begin{multline}
J_{\lambda ,\alpha }\left( m,u\right) :=e^{-2ac^{\lambda
}}\int\limits_{Q_{T}}\left[ L_{1}(u,m)^{2}+qdL_{2}(u,m)^{2}\right] \phi
_{\lambda }^{2}\ dxdt+  \label{4.1} \\
+\alpha \int_{Q_{T}}\left( u_{x}^{2}+m_{x}^{2}+u_{xx}^{2}+m_{xx}^{2}\right)
dxdt.
\end{multline}%
Recall that $L_{1}(u,m)$ and $L_{2}(u,m)$ are the operators on the left-hand
side of \eqref{3.01} and \eqref{3.02}, respectively. The only difference
between \eqref{4.1} and \eqref{3.1} is the $H^{s}$-norm regularization term
is replaced with the $L^{2}$-norms of the first and second spatial
derivatives of $u$ and $m$. This helps to simplify the numerical
implementation while maintaining solution quality. The minimization of the
functional \eqref{4.1} is carried out using the MATLAB function \texttt{%
fmincon}. For each 11-week period, the spatial domain $[-1,1]$ is
discretized into 21 grid points, while the temporal domain $[0,1]$ is
divided into 11 time steps, corresponding to the 11 weeks. The stopping
criterion for \texttt{fmincon} is set to a first-order optimality tolerance
of $10^{-5}$. We refer to MATLAB 2024a's Optimization Toolbox documentation \cite{Matlab2024a} for further details on the first-order optimality tolerance. The fixed parameters are specified as follows: 
\begin{equation*}
K(x,y)\equiv 1,\ \lambda =1,\ a=1.1,\ c=3,\ d=1,\ \alpha =10^{-4}.
\end{equation*}
{See Remark \ref{remark:kernel} for a discussion on the choice of kernel function $K(x,y)$}.

While the theory requires choosing $\lambda$ sufficiently large for the CWF, our computations use $\lambda=1$. This choice is consistent with prior convexification studies, where $\lambda\in[1,5]$ yielded stable reconstructions and accurate results; see, for example, \cite{MFGCAMWA,MFG8IPI,MFG2}. A similar practical note appears in \cite[Remark~7.2, Section~2.3.2]{Kliba2025}. In our setting, $\lambda=1$ provides reliable numerics without compromising the theoretical guarantees used to derive the algorithm.

Some numerical refinements were made to utilize the data more effectively
compared to the basic implementation studied in \cite{Kliba2025}. These
refinements include:

\begin{enumerate}
\item In addition to the zero Neumann boundary conditions and the initial
conditions 
\begin{equation*}
u(x,0)=u_{\text{est}}(x,0)\text{ and }m(x,0)=m_{\text{data}}(x,0),
\end{equation*}%
we impose the additional constraints 
\begin{equation}
m(x,h)=m_{\text{data}}(x,h)\text{ and }m(x,2h)=m_{\text{data}}(x,2h)
\label{additional_constraints}
\end{equation}%
when minimizing the functional \eqref{4.1}. In the above constraints, $u_{%
\text{est}}$ is the estimated value of $u(x,0)$, and $m_{\text{data}}$ is
the observed sentiment data. These constraints are consistent with the fact
that we have already utilized the first three weeks of data to estimate $%
u(x,0)$ as described in Section \ref{sec:initial_value_estimation}. 

\item The initial guess for \texttt{fmincon} is a pair $(u_{\text{ini}},m_{%
\text{ini}})$ defined as follows: 
\begin{gather}
u_{\text{ini}}(x,t)=u(x,0),\quad \forall t\in \lbrack 0,T], \\
m_{\text{ini}}(x,t)=
\begin{cases}
m(x,t) & \text{if }t\in \{0,h,2h\}, \\ 
\frac{1}{3}\left( m(x,0)+m(x,h)+m(x,2h)\right)  & \text{if }t>2h.%
\end{cases} \label{ini_guess}
\end{gather}
These intial guesses are already known from the data and the estimation
procedure, which does not require any further assumptions. {See Remark \ref{remark:initial_guess} for further discussion on the initial guess.}

\item Prior to the use, the raw sentiment data are smoothed using cubic
spline interpolation to reduce roughness. Minor adjustments are made near
the endpoints $x=\pm 1$ to ensure the satisfaction of the Neumann boundary
conditions. Consequently, the convexification solution for the first three
weeks do not match the raw data exactly but instead match the preprocessed
versions.
\end{enumerate}

\subsection{Results}

Through extensive numerical experiments to calibrate the parameters $\beta $%
, $r$, and $u(-1,0)$, we identified optimal values (within our computational
capacity) for each time period. These values are summarized in Table \ref%
{tab:parameters}.

\begin{table}[ht!]
\centering
\begin{tabular}{|c|c|c|c|c|}
\hline
No. & Period & $u(-1,0)$ & $\beta$ & $r$ \\ \hline
1 & Mar 2 -- May 17, 2020 & 2 & 0.25 & 50 \\ 
2 & May 18 -- Aug 2, 2020 & 1 & 0.05 & 80 \\ 
3 & Aug 3 -- Oct 18, 2020 & 2 & 0.5 & 75 \\ 
4 & Oct 19, 2020 -- Jan 3, 2021 & 3.8 & 0.3 & 80 \\ 
5 & Jan 4 -- Mar 21, 2021 & 3.7 & 1.5 & 100 \\ 
6 & Mar 22 -- Jun 6, 2021 & 5 & 3 & 200 \\ 
7 & Jun 7 -- Aug 22, 2021 & 1.6 & 2.75 & 300 \\ 
8 & Aug 23 -- Nov 7, 2021 & 2.6 & 1.5 & 125 \\ 
9 & Nov 8, 2021 -- Jan 23, 2022 & 3 & 0.75 & 75 \\ 
10 & Jan 24 -- Apr 10, 2022 & 4 & 2.5 & 175 \\ \hline
\end{tabular}%
\caption{Calibrated parameters for each period.}
\label{tab:parameters}
\end{table}

The calibrated diffusion coefficients $\beta$ exhibit substantial variation
across periods, ranging from $0.05$ to $3.0$. Similarly, the drift
coefficients $r$ vary considerably, from $50$ to $300$. This reflects
differing levels of sentiment volatility during distinct phases of the
pandemic.

\begin{remark}
The calibrated parameter values reported in Table~\ref{tab:parameters} may
serve as useful initial estimates or reference points for future studies
involving sentiment dynamics or related MFG applications. Some heuristic
approaches could be employed to find better parameter values, such as, e.g.
machine learning. However, these methods would require significant
computational resources and are beyond the scope of the current work, which
focuses on demonstrating the feasibility of MFG-based sentiment forecasting
using convexification.
\end{remark}

\begin{figure}[H]
\centering
\includegraphics[page=7, width=0.32\textwidth]{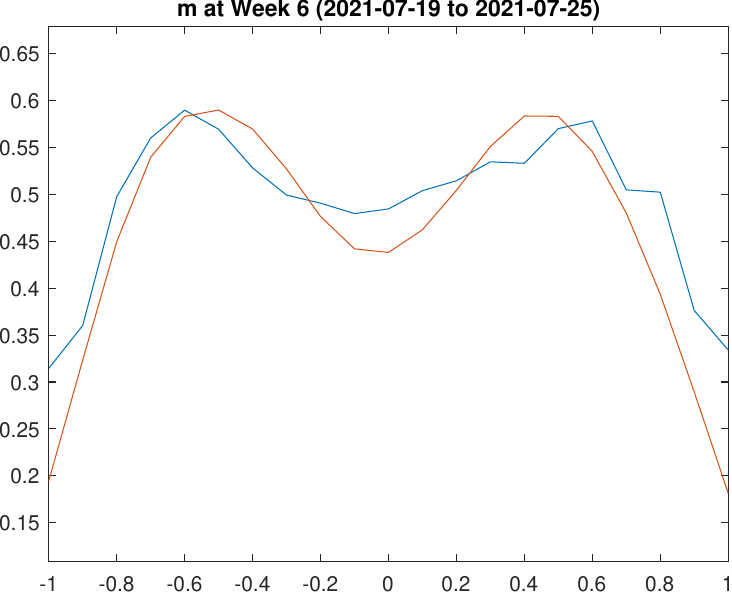}  %
\includegraphics[page=8, width=0.32\textwidth]{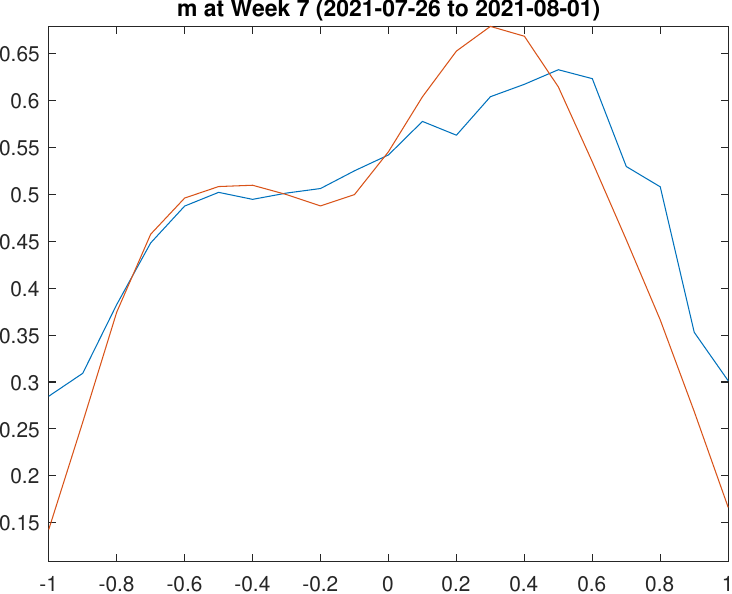}  %
\includegraphics[page=9, width=0.32\textwidth]{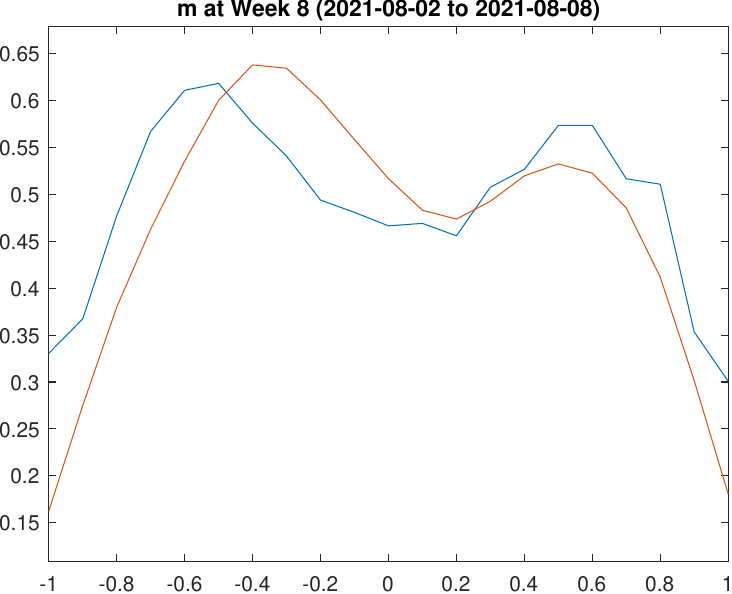}  
\caption{Convexification solution versus observed sentiment data for
weeks 6--8 in Period 7. The blue curve represents the convexification
solution, and the red curve represents the observed data. The solution captures the balanced sentiment dynamics in week 6, the surge in positive sentiment in week 7, and the surge in negative sentiment in week 8.}
\label{fig:rep1}
\end{figure}

For each of the ten periods described above, we followed the forecasting
procedure described in Sections \ref{sec:forecasting_methodology} and \ref%
{sec:Implementation} in order to compute the solution $(u,m)$ minimizing the
functional $J_{\lambda,\alpha}$ over the given period. See Figure \ref%
{fig:rep1} for a comparison of the solution $m$ with the actual values from
the dataset during Period 7 over weeks 6--8. Additionally, Figures \ref%
{fig:forecast_1m}--\ref{fig:forecast_10m} in \ref{sec:appendix} illustrate
the solution $m(x,t)$ (and the actual values from the dataset) for five
representatives of the ten 11-week periods . Note that the convexification
solution matches the preprocessed data exactly for the first three weeks of
each period due to the additional constraints \eqref{additional_constraints}%
. For subsequent weeks, the solution generally aligns well with observed
data, capturing overall trends and key features, although some discrepancies
arise due to the inherent complexity of sentiment dynamics. For example,
Figure \ref{fig:rep1} demonstrates that the solution is able to capture the
dramatic back-and-forths between positive and negative sentiment during
weeks 6--8 in Period 7.

Some other notable examples include:

\begin{enumerate}
\item Period 1 (Figure \ref{fig:forecast_1m}): The solution successfully
captures the transition from predominantly positive sentiment in week 4 to
more balanced sentiment from week 5 onward. 

\item Period 4 (Figure \ref{fig:forecast_4m}): The solution is able to
capture the surge in positive sentiment during week 3. 

\item Period 10 (Figure \ref{fig:forecast_10m}): The solution effectively
tracks the gradual increase in positive sentiment across weeks 6--10.
\end{enumerate}

\begin{figure}[H]
\centering
\includegraphics[page=23, width=0.5\textwidth]{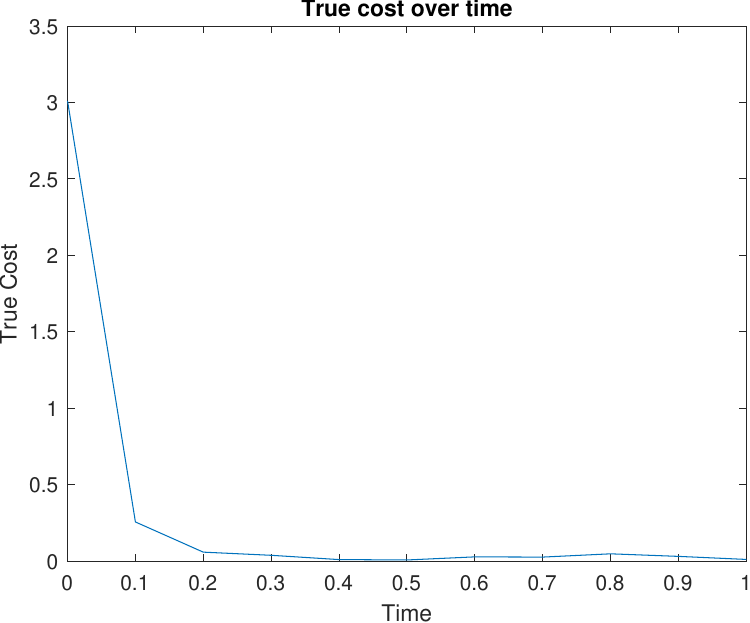}  
\caption{The true cost \eqref{true_cost} for Period 1.}
\label{fig:rep2}
\end{figure}

For each time period, we also computed the \textit{True Cost} over time, 
which measures how well the sentiment densities and value functions from
convexification satisfy the MFG system. Specifically, this cost is a
relative, unweighted, and unpenalized version of the functional $J_{\lambda
,\alpha }$ evaluated at the convexification solution $(u_{\text{sol}},m_{%
\text{sol}})$ defined as 
\begin{equation}
\text{True Cost }(t):=\left[ \frac{\int\limits_{\Omega }\left[ L_{1}(u_{%
\text{sol}},m_{\text{sol}})^{2}+L_{2}(u_{\text{sol}},m_{\text{sol}})^{2}%
\right] dx}{\int_{\Omega }\left[ u(x,0)^{2}+m(x,0)^{2}\right] dx}\right]
^{1/2}.  \label{true_cost}
\end{equation}%
Smaller values indicate better satisfaction of the governing equations. As a
representative example, Figure \ref{fig:rep2} shows this true cost value for
Period 1. Starting from $t=0.2$ (week 2), the true cost remains low,
indicating that the convexification solution well satisfies the MFG system.
We note that this behavior of the true cost is similar to that in our
previous work \cite{Kliba2025} for simulated data. See \ref%
{sec:appendix} for the true cost figures for each of our representative
examples.

Since the dataset does not include actual value functions $u(x,t)$ for
validation, we only present the predicted value function by convexification
for one representative period in Figure \ref{fig:forecast_1u} in \ref%
{sec:appendixB}. The validity of the estimated value functions is supported
by consistently low true cost values throughout the forecast horizon.

\begin{figure}[H]
\centering
\includegraphics[page=1, width=0.5\textwidth]{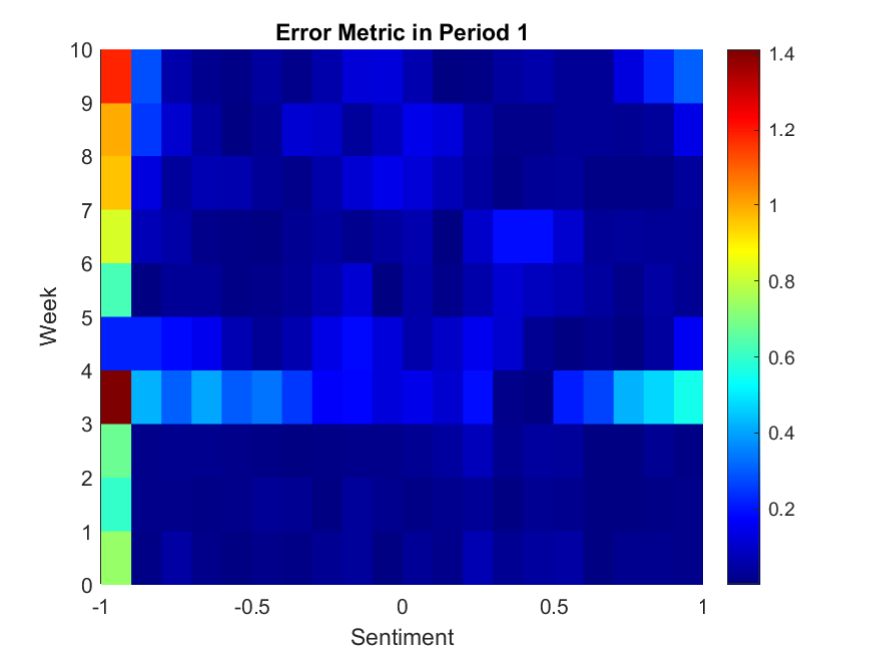}  
\caption{Error metric \eqref{error_metric} for Period 1. Each tile's color
represents the metric value for the corresponding week and sentiment. Dark
blue tiles are where the metric is smallest, while dark red tiles are where
the metric is largest.}
\label{fig:rep3}
\end{figure}

For each of the time periods, we also computed a metric 
for the relative errors between the convexification solution and the
observed sentiment densities at each $(x,t)$. 
This relative error metric is computed as follows: 
\begin{equation}
\text{Error }(x,t):=\frac{\left\vert m_{\text{sol}}(x,t)-m_{\text{data}%
}(x,t)\right\vert }{\left\vert m_{\text{data}}(x,t)\right\vert }.
\label{error_metric}
\end{equation}%
Across all ten periods, this error metric value remains below $25\%$ for
most $(x,t)$ pairs. Larger errors typically appear near the boundary points $%
x=\pm 1$, where the observed sentiment densities deviate from the Neumann
boundary conditions. For instance, Figure \ref{fig:rep3} shows the error
metric for Period 1. See Figure \ref{fig:errors} in \ref{sec:appendixC} for
all time periods.

\begin{remark}
Throughout each period, the true cost remains consistently low, indicating
that the convexification solution well satisfies the MFG system. Combined
with the relatively small error between the convexification solution and
observed densities, this demonstrates the ability of the
convexification-based MFG framework to capture the underlying dynamics of
public sentiment in this dataset.
\end{remark}

\begin{figure}[H]
\centering
\includegraphics[page=5, width=0.32\textwidth]{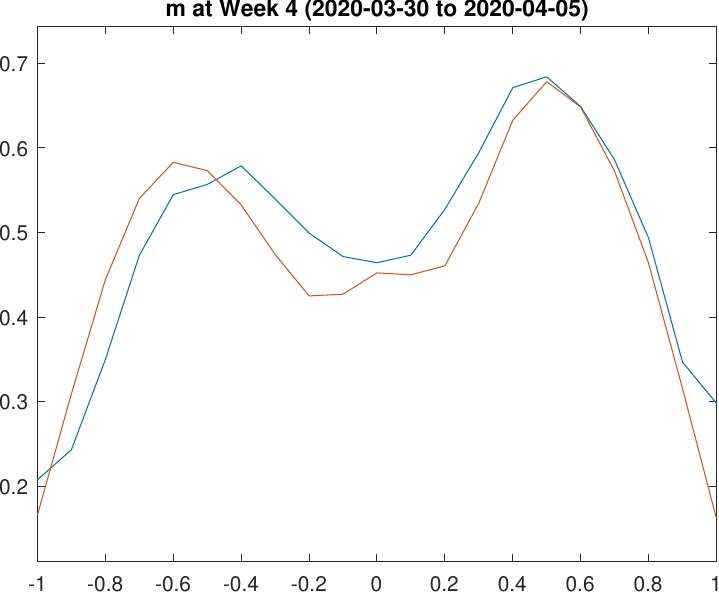} 
\includegraphics[page=6, width=0.32\textwidth]{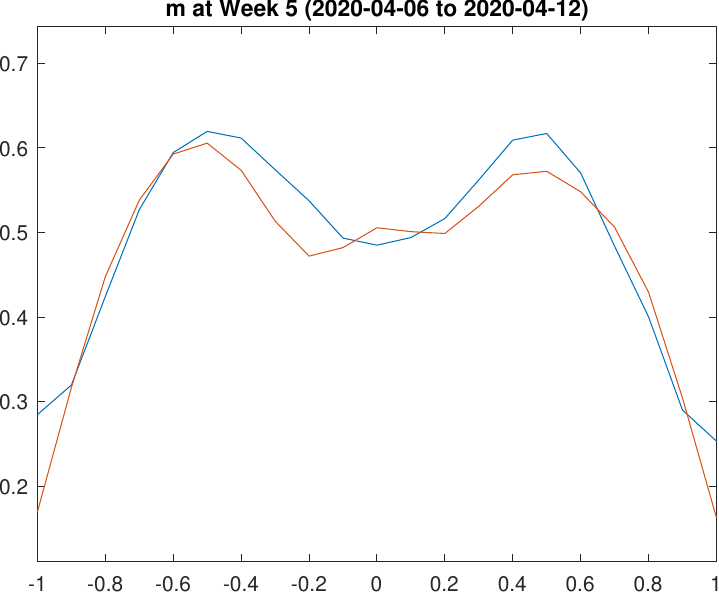}
\includegraphics[page=7, width=0.32\textwidth]{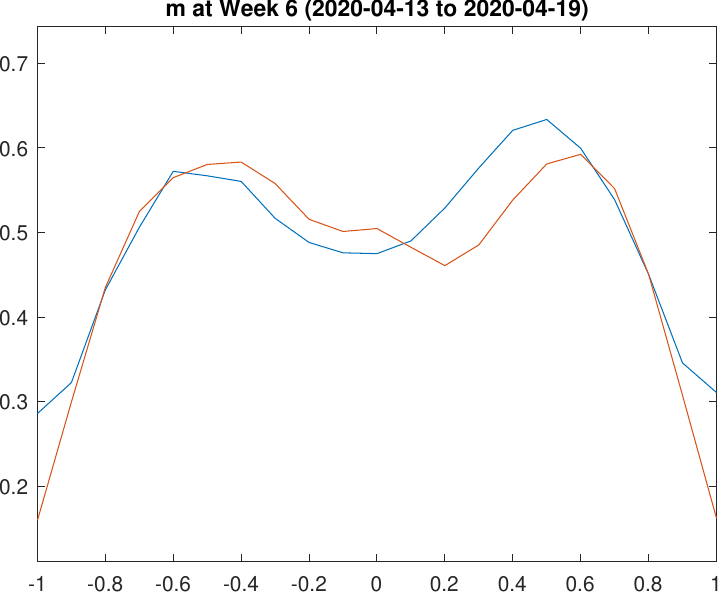} \\ \vspace{5mm}
\includegraphics[page=5, width=0.32\textwidth]{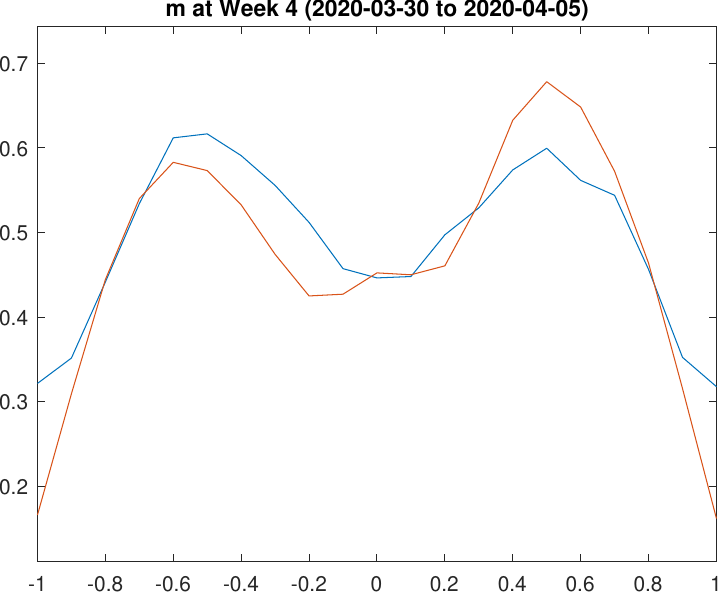}
\includegraphics[page=6, width=0.32\textwidth]{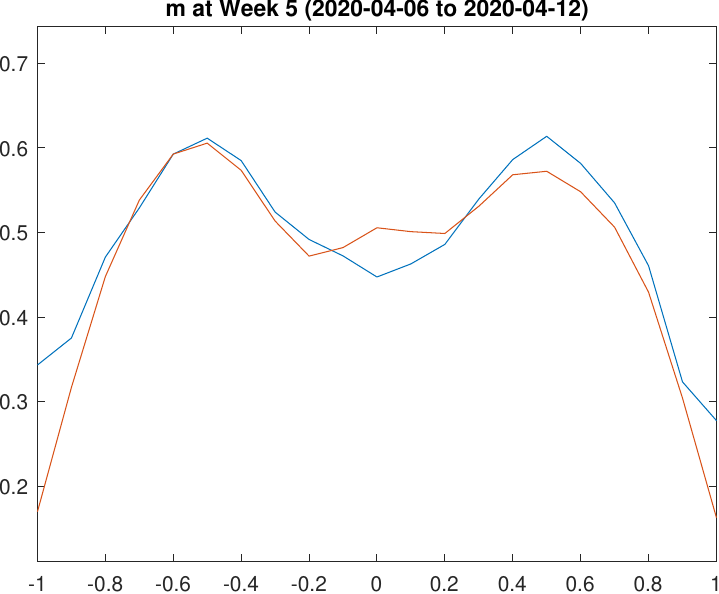} 
\includegraphics[page=7, width=0.32\textwidth]{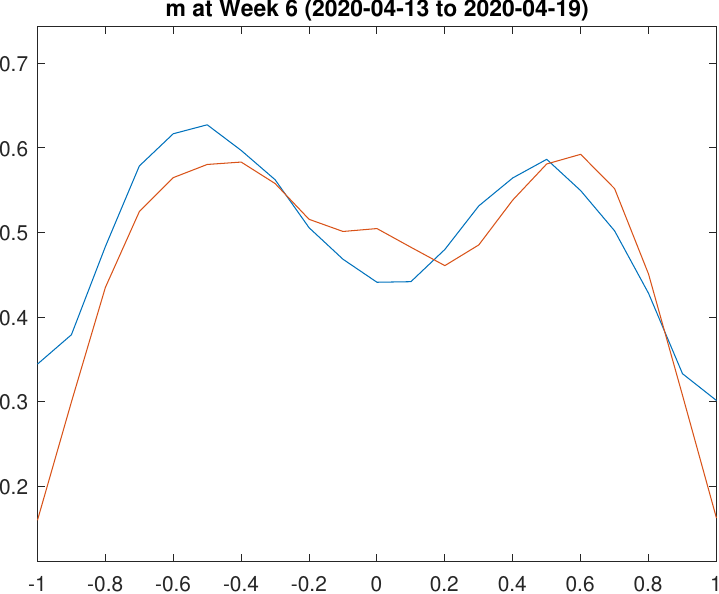}
\caption{{Convexification solutions for weeks 4--6 in Period 1 with different kernel choices versus observed data. The blue curve represents the convexification solution, and the red curve represents the observed data. The top row uses the constant kernel, while the bottom row uses the Gaussian kernel. The different kernel choices lead to different trends in the solution, with the constant kernel producing a closer fit to the observed data under the calibrated parameters.}}
\label{fig:kernel}
\end{figure}

{\begin{remark} \label{remark:kernel}
    In addition to the constant kernel $K(x,y)\equiv 1$, we also experimented with a Gaussian kernel 
    $$
        K(x,y)=\frac{1}{\sqrt{2\pi\sigma^2}}\exp \left(-\frac{(x-y)^2}{2\sigma^2} \right)
    $$ 
    with $\sigma=0.1$. Figure \ref{fig:kernel} compares the convexification solutions with these two kernel choices for weeks 4--6 in Period 1. The solution with the constant kernel (top row) shows a predomiantly positive sentiment in week 4, a balanced trend in week 5, and a slightly positive sentiment in week 6. On the other hand, the solution with the Gaussian kernel (bottom row) shows a slightly negative sentiment in weeks 4 and 6, and a balanced trend in week 5. This suggests that the choice of kernel can influence the solution's behavior. We note that, to simplify the presentation, the choice of the other parameters remains the same for both kernel choices (see their values Section \ref{sec:Implementation} and Table \ref{tab:parameters}). Under these parameters, the constant kernel aligns better with the observed data. We conjecture, however, that if our choice of parameters were adapted to the Gaussian kernel, the solution with the Gaussian kernel could also fit the observed data well. This highlights the importance of parameter calibration in conjunction with kernel choice.
\end{remark}
}

\begin{figure}[H]
\centering
\includegraphics[page=4, width=0.32\textwidth]{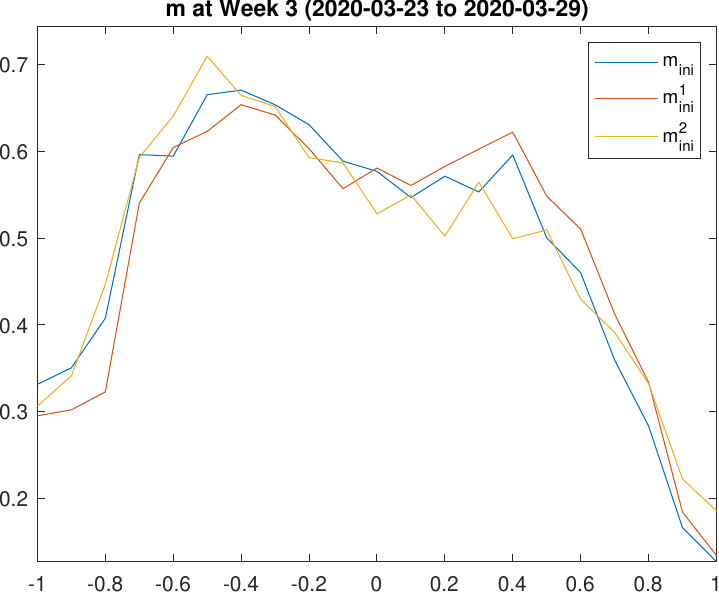}  %
\includegraphics[page=5, width=0.32\textwidth]{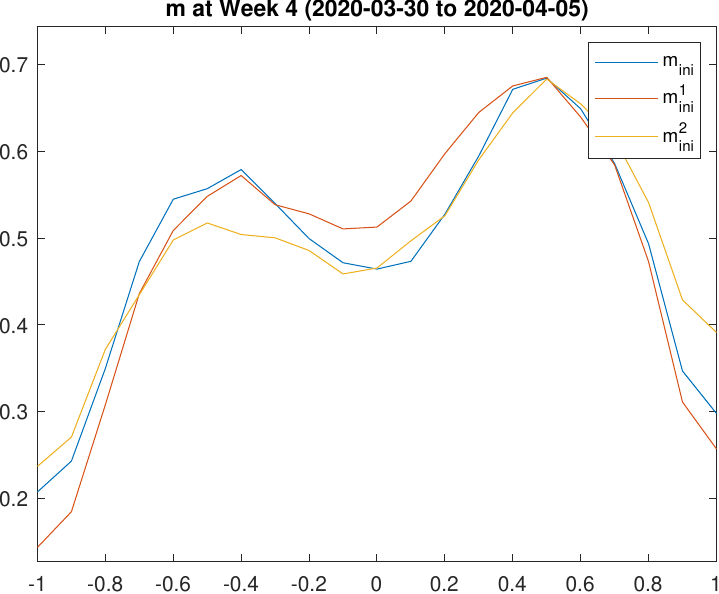}  %
\includegraphics[page=6, width=0.32\textwidth]{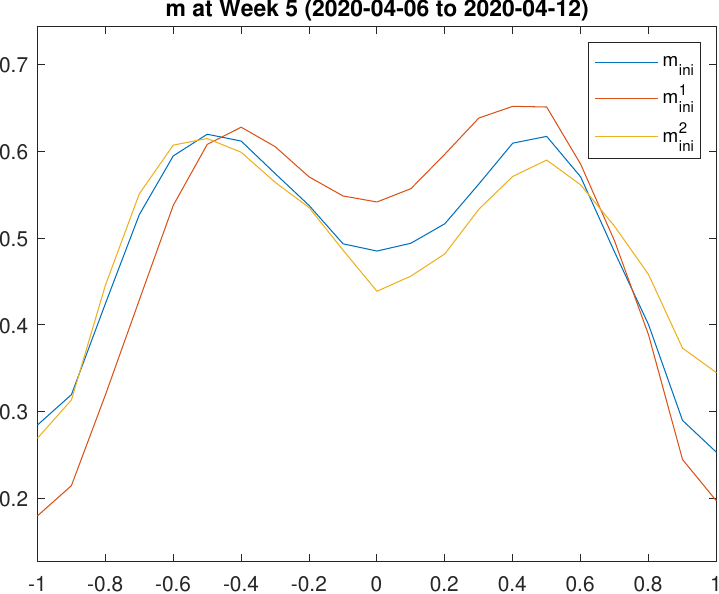}  
\caption{{Convexification solutions for weeks 3--5 in Period 1 with three different initial guesses: $m_{\text{ini}}$, $m^1_{\text{ini}}$, and $m^2_{\text{ini}}$. The blue curve represents the solution with $m_{\text{ini}}$, the orange curve represents the solution with $m^1_{\text{ini}}$, and the yellow curve represents the solution with $m^2_{\text{ini}}$. The solutions show very close overall trends and key features, such as a predominantly positive trend in week 3, a predominantly negative trend in week 4, and the return to a balanced trend in week 5.}}
\label{fig:guesses}
\end{figure}

{\begin{remark} \label{remark:initial_guess}
    Convexification guarantees the existence of a unique minimizer regardless of the initial guess used in the optimization. However, from the perspective of numerical optimization, the choice of the initial guess can influence the convergence speed and the quality of the solution obtained within a reasonable computational time. In our implementation, we want to best utilize the available data to construct a well-informed initial guess. This approach helps to guide the optimization process towards a more accurate solution while ensuring that it remains computationally efficient. Aside from the initial guess defined in \eqref{ini_guess}, we also experimented with other initial guesses for $m(x,t)$, such as 
    $$
        m^1_{\text{ini}}(x,t)=
        \begin{cases}
        m(x,t) & \text{if }t\in \{0,h,2h\}, \\ 
        \frac{1}{2}\left( m(x,0)+m(x,2h)\right)  & \text{if }t>2h,
        \end{cases}
    $$
    and
    $$
        m^2_{\text{ini}}(x,t)=
        \begin{cases}
        m(x,t) & \text{if }t\in \{0,h,2h\}, \\ 
        m(x,h)  & \text{if }t>2h.
        \end{cases}
    $$  
    Figure \ref{fig:guesses} compares convexification solutions from these alternative guesses against the original $m_{\text{ini}}$ defined in \eqref{ini_guess} for a representative period. While minor differences exist, they do not contradict theory: at an optimality tolerance of $10^{-5}$ (see Section \ref{sec:Implementation} for details), multiple nearly-optimal solutions may satisfy the stopping criterion. Importantly, all solutions exhibit consistent overall trends and key features, demonstrating the robustness of convexification to the choice of initial guess.
\end{remark}
}

\section{Summary and Future Directions}

\label{sec:discussion_future_directions}

\subsection{Summary of Findings}

This work demonstrates the practical potential of MFG models combined with
convexification-based numerical methods for forecasting public sentiment
dynamics. Using real-world sentiment data from social media responses to CDC
tweets during the COVID-19 pandemic (March 2020--April 2022), we have demonstrated
that:

\begin{enumerate}
\item The MFG framework can reproduce key features of observed sentiment
evolution, including sudden polarity shifts, gradual trends, and changes in
distributional characteristics {(see Figures \ref{fig:rep1} and \ref{fig:forecast_1m}~-~\ref{fig:forecast_10m})}. 

\item Solutions obtained via the convexification method not only match
observed data well {(see Figures \ref{fig:rep3} and \ref{fig:errors})} but also satisfy the governing MFG equations with small
residuals, as evidenced by the low true cost values {(see Figure \ref{fig:rep2} and the true cost values in Figures \ref{fig:forecast_1m}~-~\ref{fig:forecast_10m})}.
\end{enumerate}

{These results provide proof-of-concept evidence that appropriately parameterized MFG models can capture complex temporal patterns in public sentiment -- a statement that was previously a theoretical hypothesis but now has become empirically verified with real-life data. The convexification method offers both theoretical rigor and numerical efficiency as a practical approach for solving the MFG system. As the HJB equation with initial condition for $u$ is not stable, the convexification method helps to maintain a balance between stability and accuracy for $t \in [0,1]$. Moreover, our experiments with different interaction kernels and initial guesses demonstrate robustness of convexification against the choices of these quantities.}

The main limitation of the current study is that the coefficients $\beta $, $
r$, and the number  $u(-1,0)$ were determined through a manual calibration
rather than a systematic inverse problem techniques. While this approach
suffices for a proof of concept, it limits predictive capability since
values of these coefficients are unknown \textit{a priori}. Another
limitation is that the assumption of piecewise-constant coefficients over
11-week intervals may be restrictive given the dynamic nature of public
sentiment.

\subsection{Future Work}

To overcome the above limitations and advance MFG-based sentiment
forecasting, two key directions for future research are identified:

\begin{enumerate}
\item \textit{Systematic coefficient identification:} Developing rigorous
inverse problem solvers to determine $\beta$ and $r$ directly
from measured or historical data is a priority. This will enable genuine
forecasting. {The problem of determining these coefficients is an extremely challenging coefficient inverse problem due to its non-linearity, ill-posedness, and non-convexity when tackled with a standard optimization approach}.  The convexification framework has been successfully applied to
a variety of coefficient inverse problems in the MFG theory \cite{MFGbook,
MFGCAMWA, MFG8IPI}, suggesting feasibility. 

\item \textit{Multi-dimensional extensions:} Extending to higher-dimensional
sentiment spaces to provide richer representations of public sentiment.

\item {\textit{Interaction kernel identification:} Develop a principled approach for inferring an interaction kernel that best reflects the underlying dynamics in the dataset. In the present work, we use a constant kernel for simplicity, but a systematic methodology for kernel identification would allow the model to capture more complex interaction structures and improve interpretability.}
\end{enumerate}

\subsection{Concluding Remarks}

This study provides strong evidence that MFG models offer a principled and
effective framework for modeling and forecasting public sentiment dynamics.
The convexification method ensures that stable solutions to the MFG system
can be found even with noisy or incomplete data. With continued development
of systematic parameter identification techniques, MFG-based forecasting
with convexification has the potential to become a valuable tool for
understanding and predicting behavior in social systems.

\section*{Acknowledgment}

This research is supported by the U.S. National Science Foundation Grant DMS~2436227.

\appendix

\section{Convexification Solution for Representative Periods}

\label{sec:appendix}

\begin{figure}[H]
\centering
\includegraphics[page=1, width=0.32\textwidth]{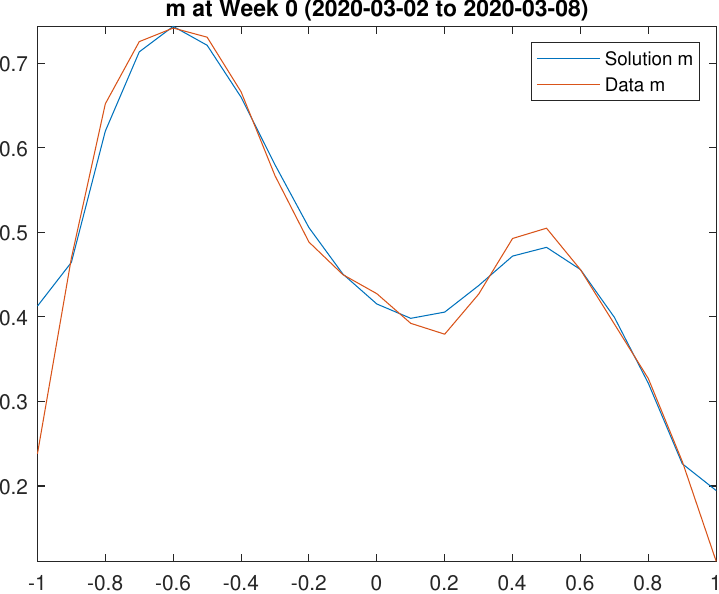} 
\includegraphics[page=2, width=0.32\textwidth]{results_period_1.pdf}
\includegraphics[page=3, width=0.32\textwidth]{results_period_1.pdf} \\ \vspace{5mm}
\includegraphics[page=4, width=0.32\textwidth]{results_period_1.pdf}
\includegraphics[page=5, width=0.32\textwidth]{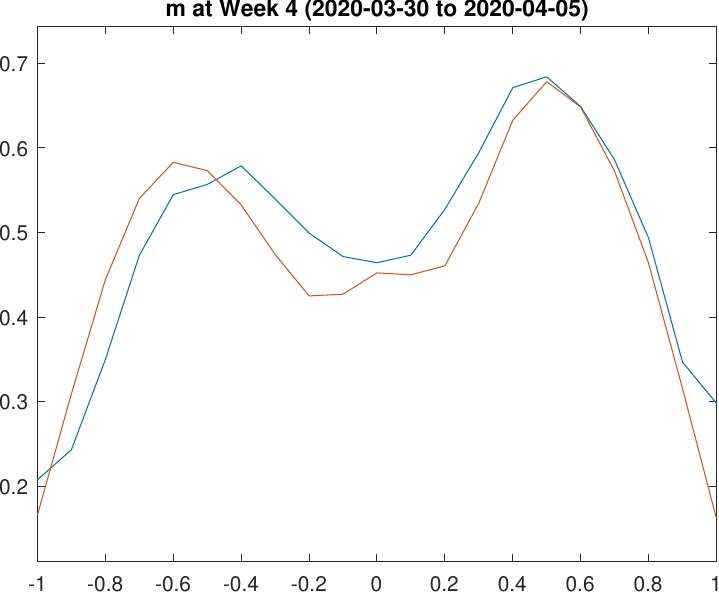} 
\includegraphics[page=6, width=0.32\textwidth]{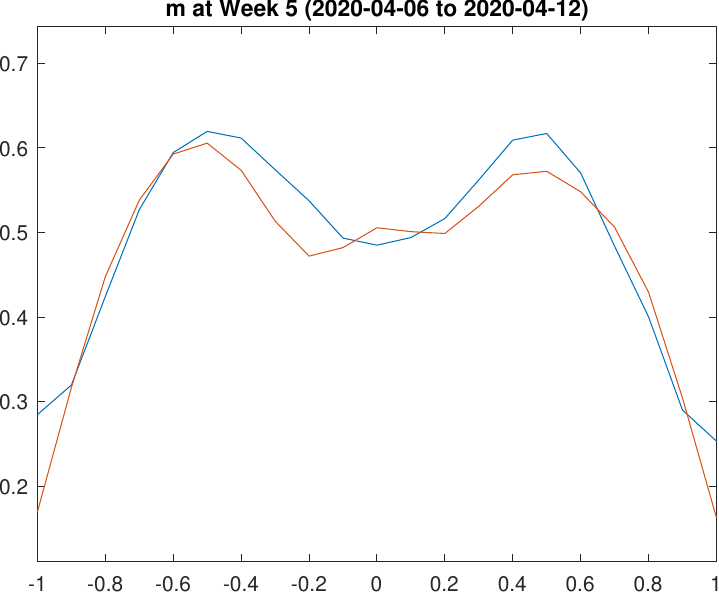} \\ \vspace{5mm}
\includegraphics[page=7, width=0.32\textwidth]{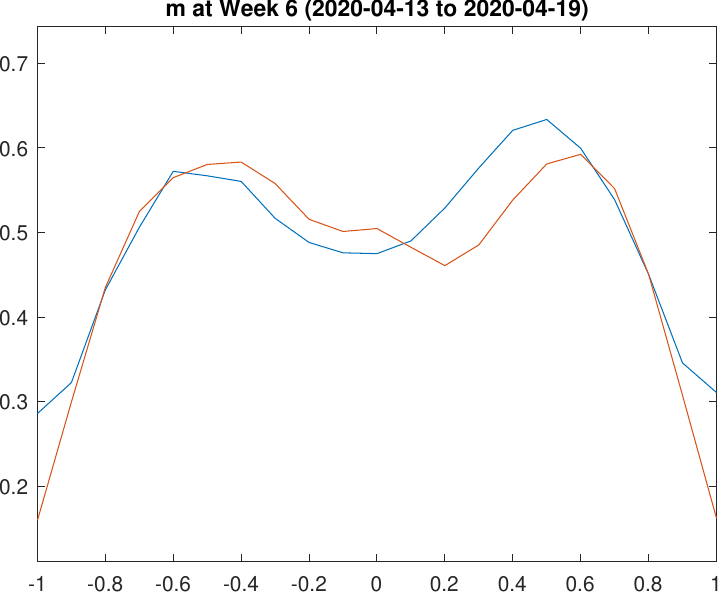}
\includegraphics[page=8, width=0.32\textwidth]{results_period_1.pdf}
\includegraphics[page=9, width=0.32\textwidth]{results_period_1.pdf} \\ \vspace{5mm}
\includegraphics[page=10, width=0.32\textwidth]{results_period_1.pdf}
\includegraphics[page=11, width=0.32\textwidth]{results_period_1.pdf}
\includegraphics[page=23, width=0.32\textwidth]{results_period_1.pdf}
\caption{Period 1: March 2 to May 17, 2020.}
\label{fig:forecast_1m}
\end{figure}

\begin{figure}[p]
\centering
\includegraphics[page=1, width=0.32\textwidth]{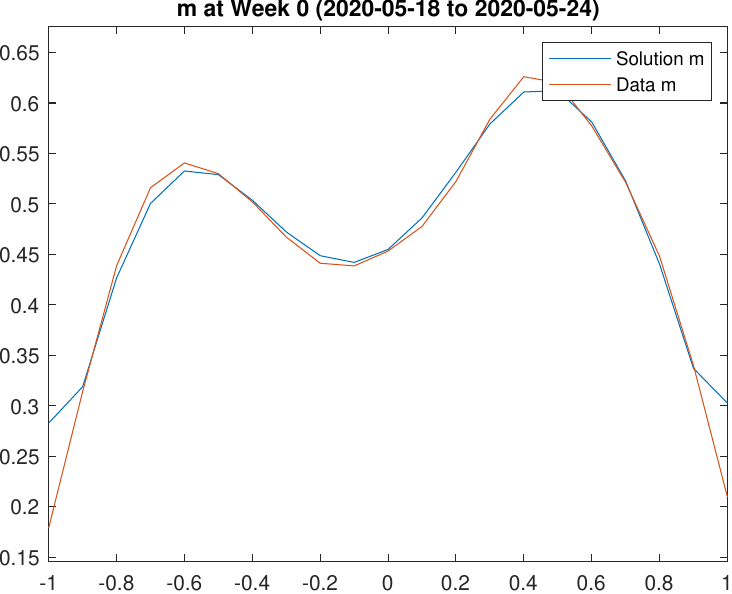} 
\includegraphics[page=2, width=0.32\textwidth]{results_period_2.pdf}
\includegraphics[page=3, width=0.32\textwidth]{results_period_2.pdf} \\ \vspace{5mm}
\includegraphics[page=4, width=0.32\textwidth]{results_period_2.pdf}
\includegraphics[page=5, width=0.32\textwidth]{results_period_2.pdf} 
\includegraphics[page=6, width=0.32\textwidth]{results_period_2.pdf} \\ \vspace{5mm}
\includegraphics[page=7, width=0.32\textwidth]{results_period_2.pdf}
\includegraphics[page=8, width=0.32\textwidth]{results_period_2.pdf}
\includegraphics[page=9, width=0.32\textwidth]{results_period_2.pdf} \\ \vspace{5mm}
\includegraphics[page=10, width=0.32\textwidth]{results_period_2.pdf}
\includegraphics[page=11, width=0.32\textwidth]{results_period_2.pdf}
\includegraphics[page=23, width=0.32\textwidth]{results_period_2.pdf}
\caption{Period 2: May 18 to August 2, 2020.}
\label{fig:forecast_2m}
\end{figure}


\begin{figure}[p]
\centering
\includegraphics[page=1, width=0.32\textwidth]{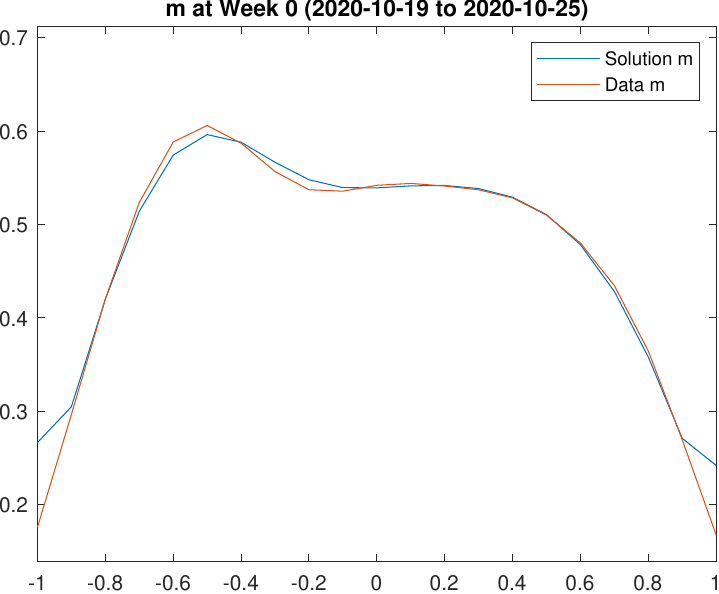} 
\includegraphics[page=2, width=0.32\textwidth]{results_period_4.pdf}
\includegraphics[page=3, width=0.32\textwidth]{results_period_4.pdf} \\ \vspace{5mm}
\includegraphics[page=4, width=0.32\textwidth]{results_period_4.pdf}
\includegraphics[page=5, width=0.32\textwidth]{results_period_4.pdf} 
\includegraphics[page=6, width=0.32\textwidth]{results_period_4.pdf} \\ \vspace{5mm}
\includegraphics[page=7, width=0.32\textwidth]{results_period_4.pdf}
\includegraphics[page=8, width=0.32\textwidth]{results_period_4.pdf}
\includegraphics[page=9, width=0.32\textwidth]{results_period_4.pdf} \\ \vspace{5mm}
\includegraphics[page=10, width=0.32\textwidth]{results_period_4.pdf}
\includegraphics[page=11, width=0.32\textwidth]{results_period_4.pdf}
\includegraphics[page=23, width=0.32\textwidth]{results_period_4.pdf}
\caption{Period 4: October 19, 2020 to January 3, 2021.}
\label{fig:forecast_4m}
\end{figure}

\begin{figure}[p]
\centering
\includegraphics[page=1, width=0.32\textwidth]{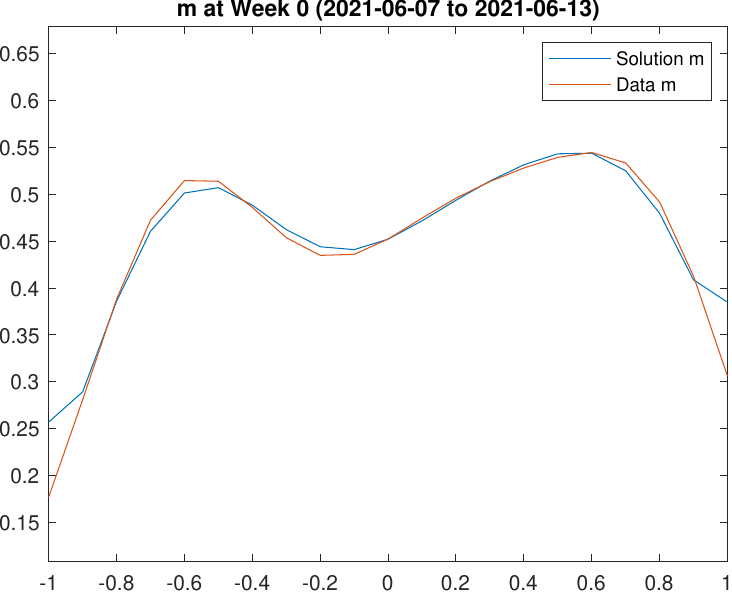} 
\includegraphics[page=2, width=0.32\textwidth]{results_period_7.pdf}
\includegraphics[page=3, width=0.32\textwidth]{results_period_7.pdf} \\ \vspace{5mm}
\includegraphics[page=4, width=0.32\textwidth]{results_period_7.pdf}
\includegraphics[page=5, width=0.32\textwidth]{results_period_7.pdf} 
\includegraphics[page=6, width=0.32\textwidth]{results_period_7.pdf} \\ \vspace{5mm}
\includegraphics[page=7, width=0.32\textwidth]{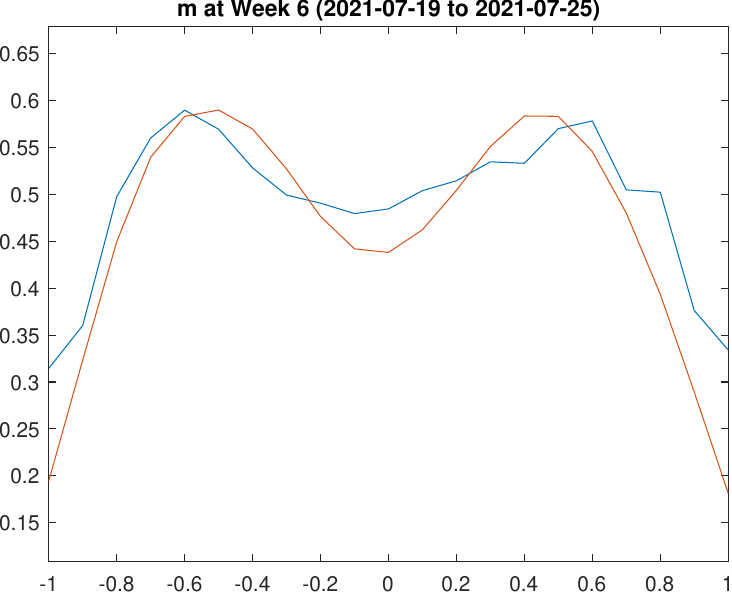}
\includegraphics[page=8, width=0.32\textwidth]{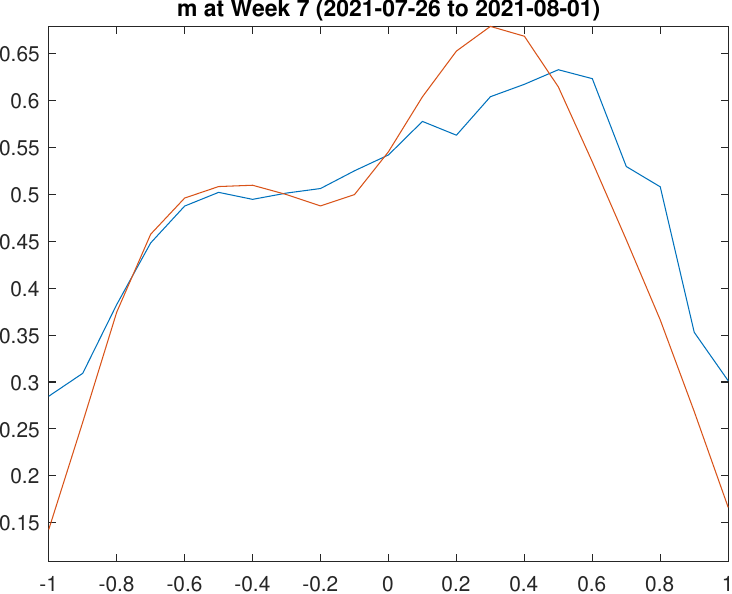}
\includegraphics[page=9, width=0.32\textwidth]{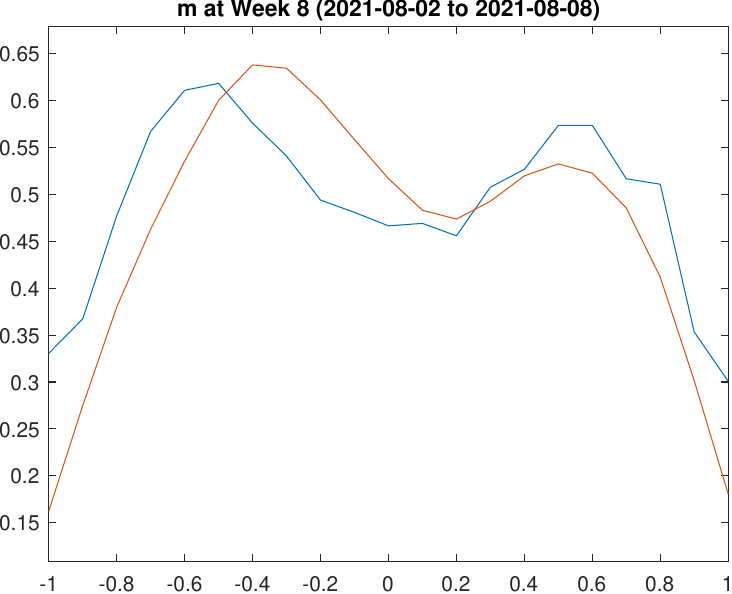} \\ \vspace{5mm}
\includegraphics[page=10, width=0.32\textwidth]{results_period_7.pdf}
\includegraphics[page=11, width=0.32\textwidth]{results_period_7.pdf}
\includegraphics[page=23, width=0.32\textwidth]{results_period_7.pdf}
\caption{Period 7: June 7 to August 22, 2021.}
\label{fig:forecast_7m}
\end{figure}

\begin{figure}[p]
\centering
\includegraphics[page=1, width=0.32\textwidth]{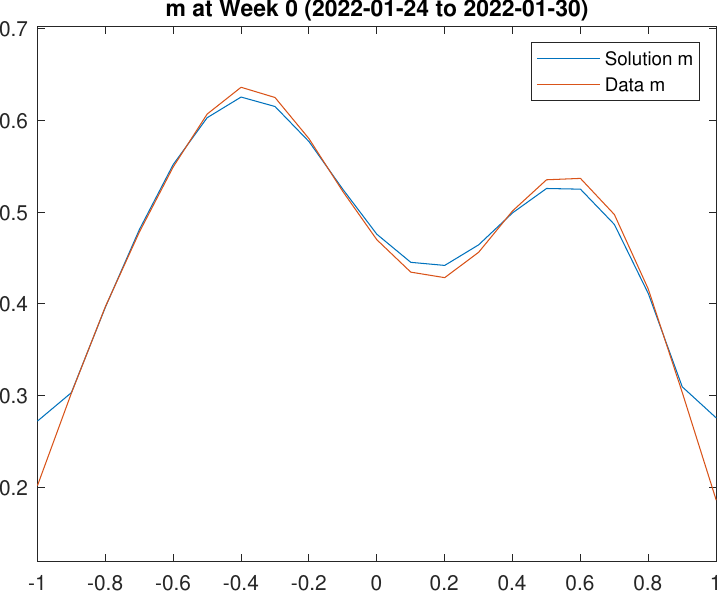} 
\includegraphics[page=2, width=0.32\textwidth]{results_period_10.pdf}
\includegraphics[page=3, width=0.32\textwidth]{results_period_10.pdf} \\ \vspace{5mm}
\includegraphics[page=4, width=0.32\textwidth]{results_period_10.pdf}
\includegraphics[page=5, width=0.32\textwidth]{results_period_10.pdf} 
\includegraphics[page=6, width=0.32\textwidth]{results_period_10.pdf} \\ \vspace{5mm}
\includegraphics[page=7, width=0.32\textwidth]{results_period_10.pdf}
\includegraphics[page=8, width=0.32\textwidth]{results_period_10.pdf}
\includegraphics[page=9, width=0.32\textwidth]{results_period_10.pdf} \\ \vspace{5mm}
\includegraphics[page=10, width=0.32\textwidth]{results_period_10.pdf}
\includegraphics[page=11, width=0.32\textwidth]{results_period_10.pdf}
\includegraphics[page=23, width=0.32\textwidth]{results_period_10.pdf}
\caption{Period 10: January 24 to April 10, 2022.}
\label{fig:forecast_10m}
\end{figure}

\section{Predicted Value Function for Period 1}
\label{sec:appendixB}

\begin{figure}[H]
\centering
\includegraphics[page=12, width=0.32\textwidth]{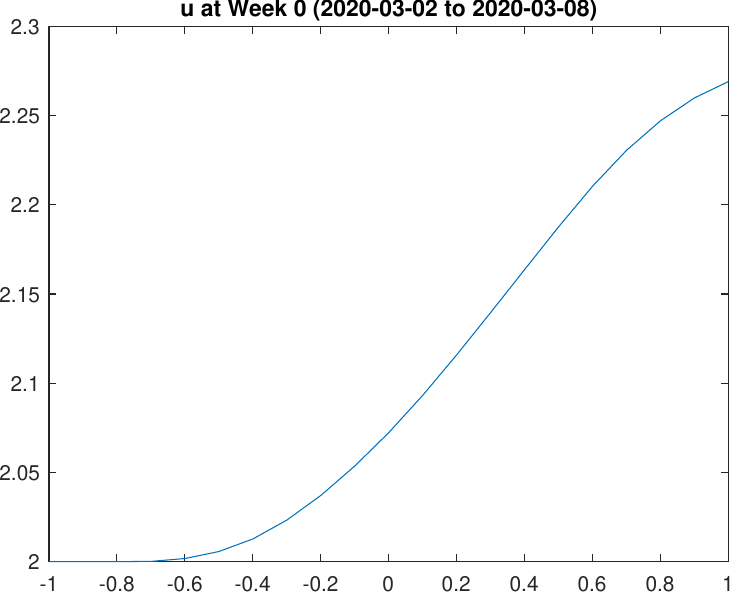} 
\includegraphics[page=13, width=0.32\textwidth]{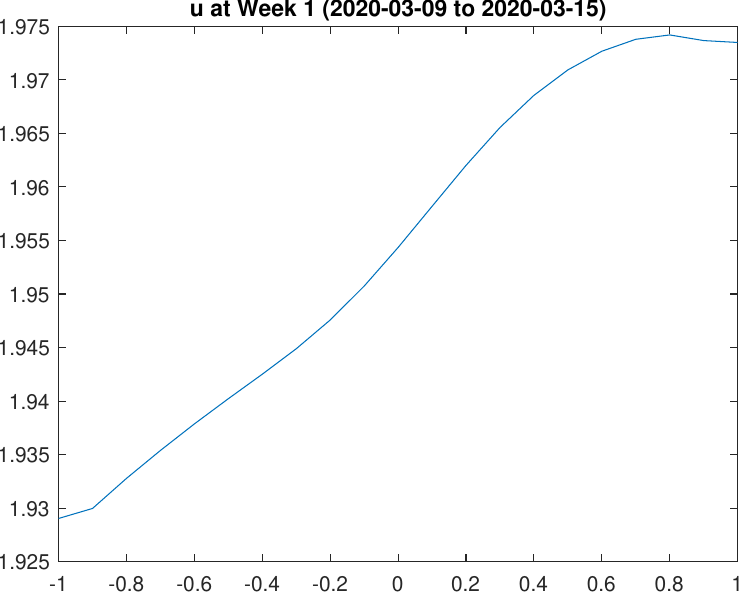}
\includegraphics[page=14, width=0.32\textwidth]{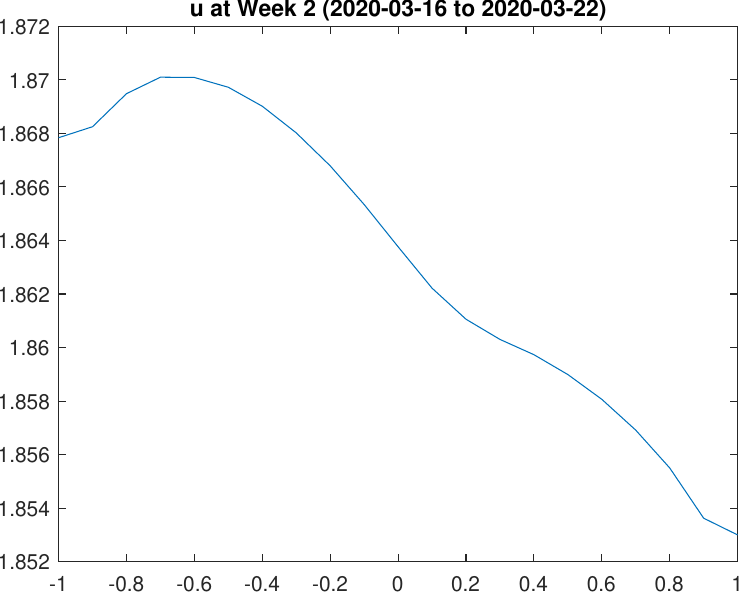} \\ \vspace{5mm}
\includegraphics[page=15, width=0.32\textwidth]{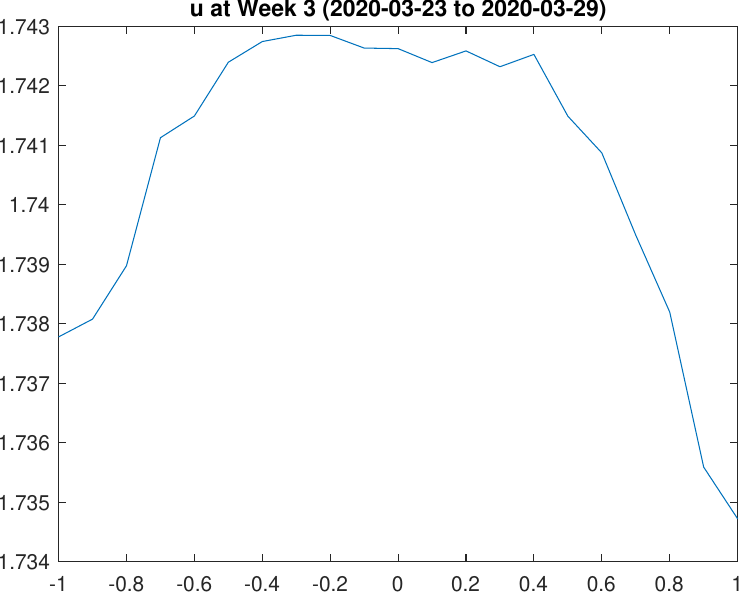}
\includegraphics[page=16, width=0.32\textwidth]{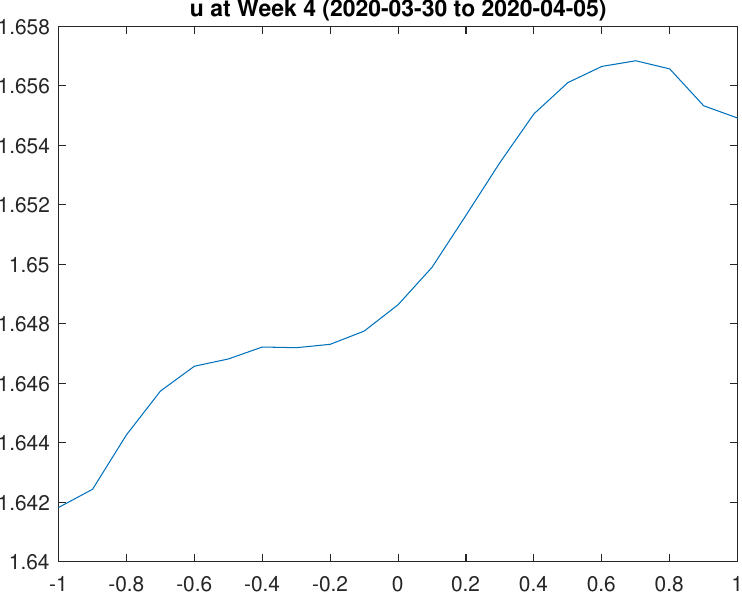} 
\includegraphics[page=17, width=0.32\textwidth]{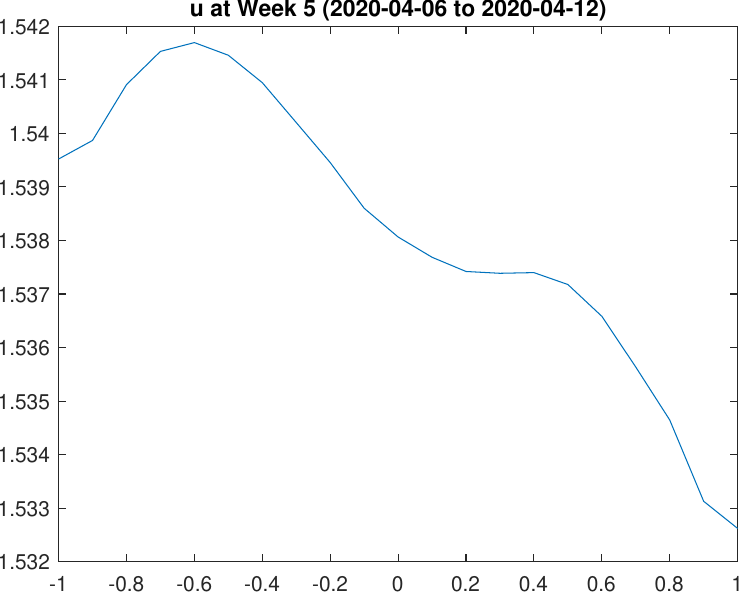} \\ \vspace{5mm}
\includegraphics[page=18, width=0.32\textwidth]{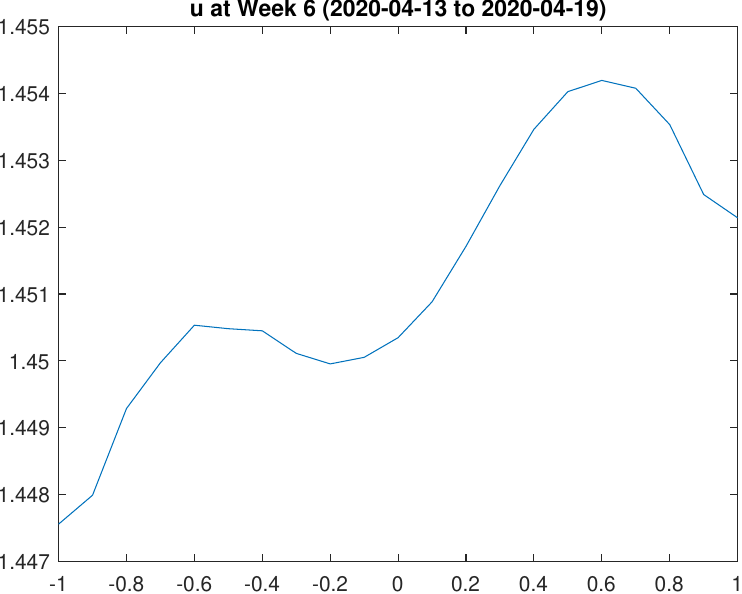}
\includegraphics[page=19, width=0.32\textwidth]{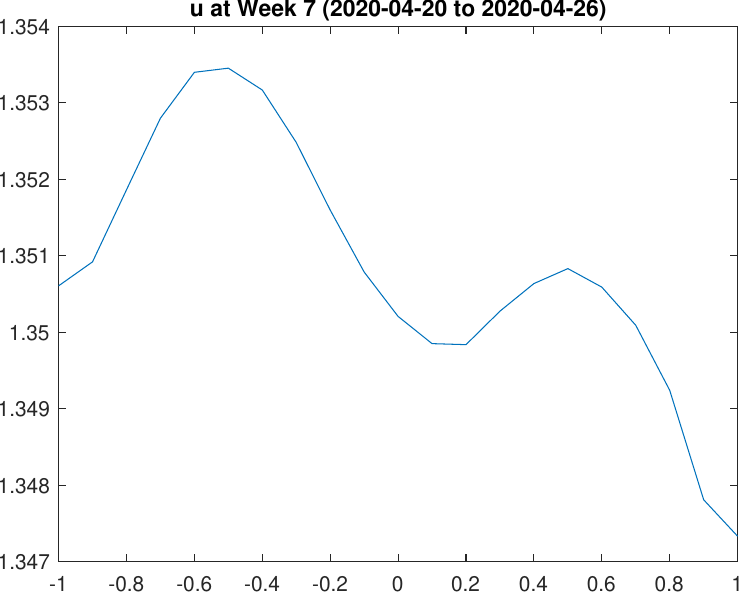}
\includegraphics[page=20, width=0.32\textwidth]{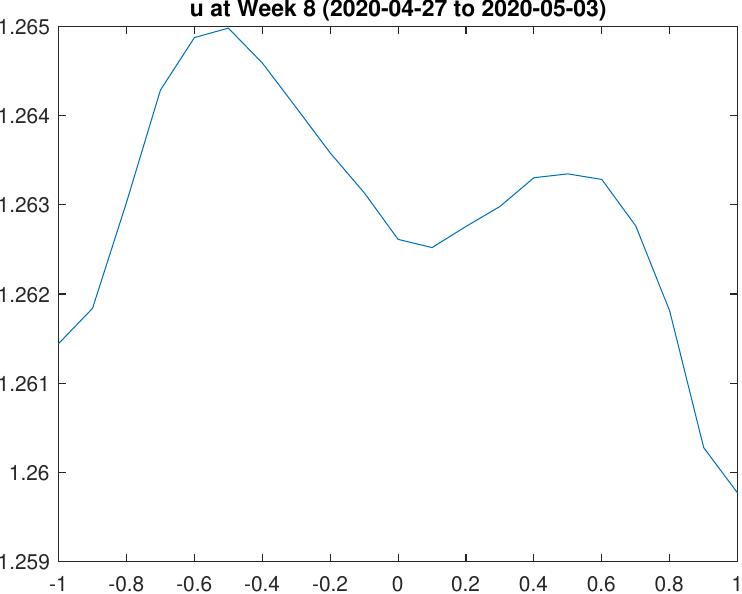} \\ \vspace{5mm}
\includegraphics[page=21, width=0.32\textwidth]{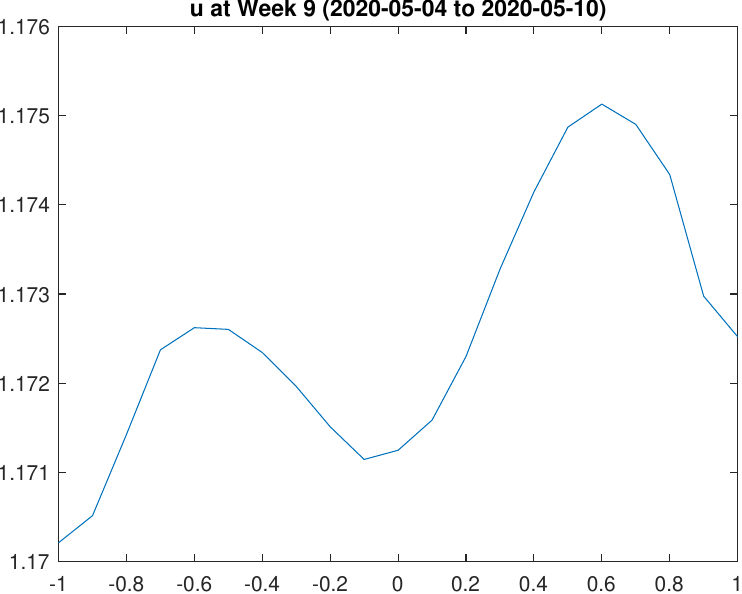}
\includegraphics[page=22, width=0.32\textwidth]{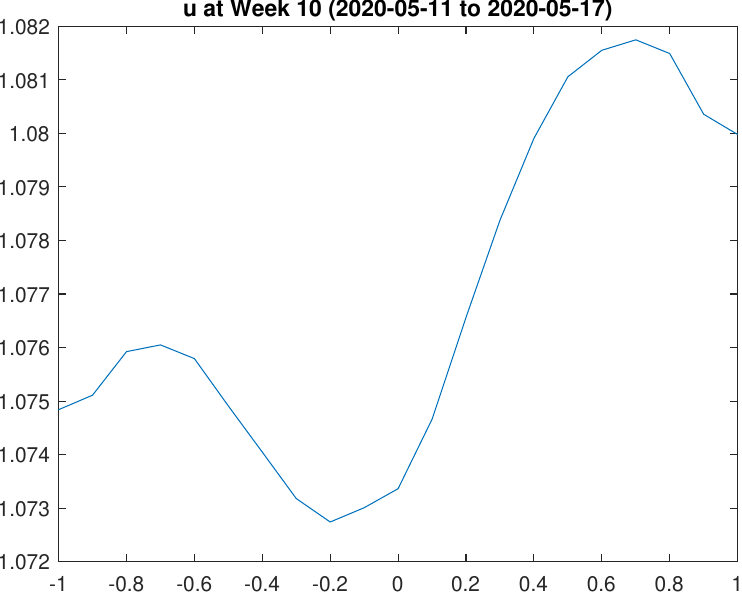}
\caption{Predicted value function $u(x,t)$ for each of the 11 weeks from March 2 to May 17, 2020.}
\label{fig:forecast_1u}
\end{figure}

\section{Overall Error Metric}
\label{sec:appendixC}

\begin{figure}[H]
\centering
\includegraphics[page=1, width=0.325\textwidth]{all_error.pdf} 
\includegraphics[page=2, width=0.325\textwidth]{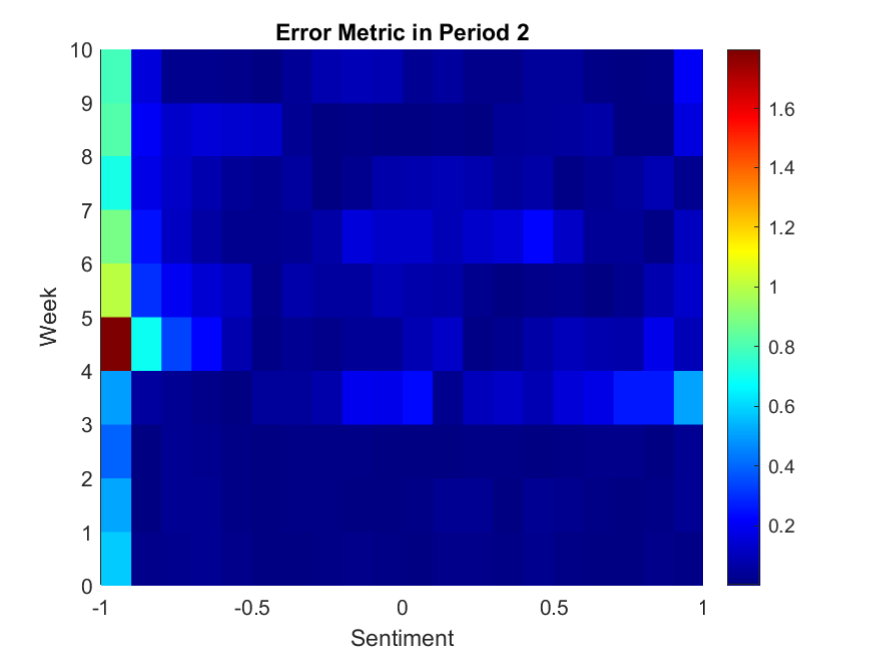} 
\includegraphics[page=3, width=0.325\textwidth]{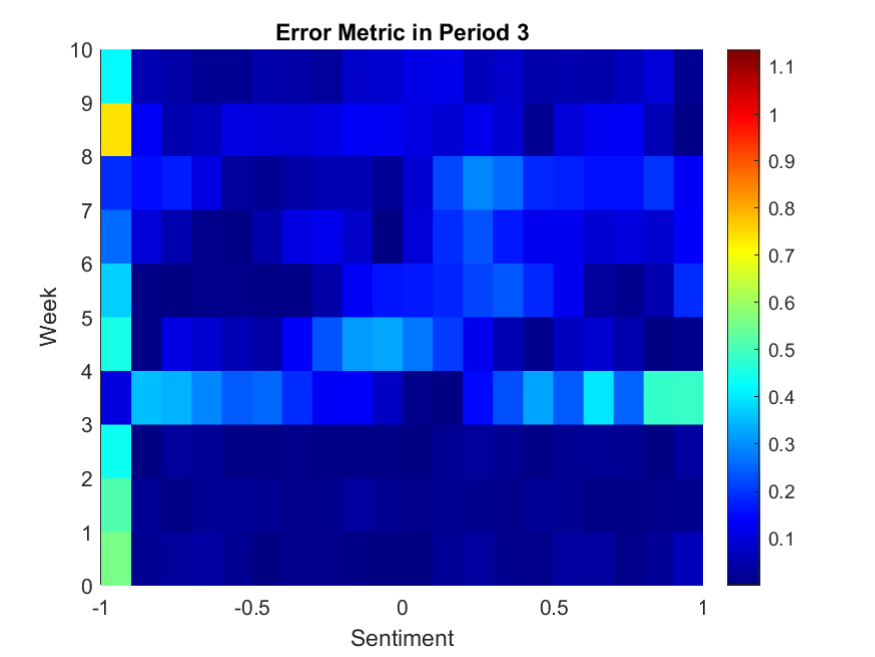} \\ \vspace{5mm}
\includegraphics[page=4, width=0.325\textwidth]{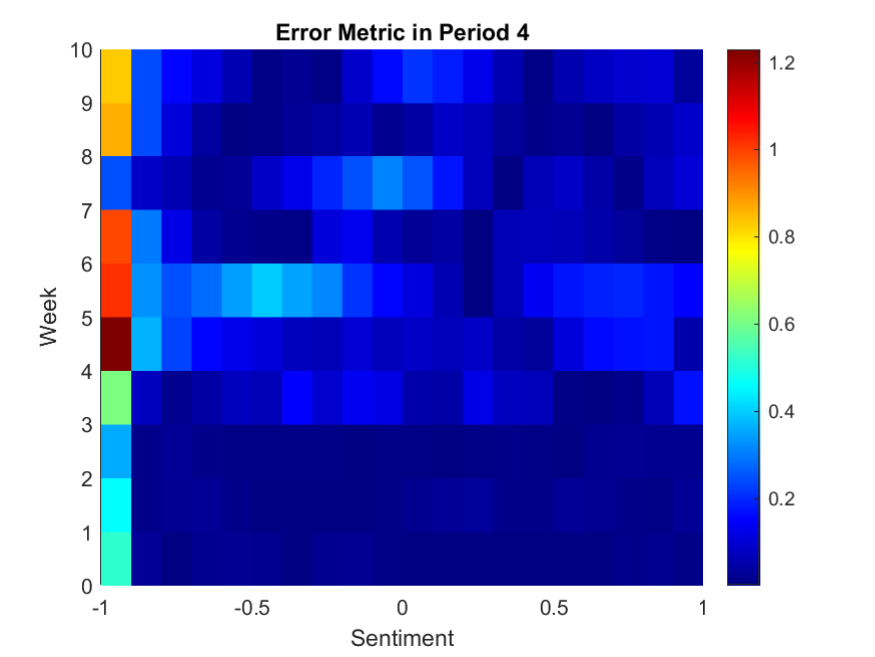} 
\includegraphics[page=5, width=0.325\textwidth]{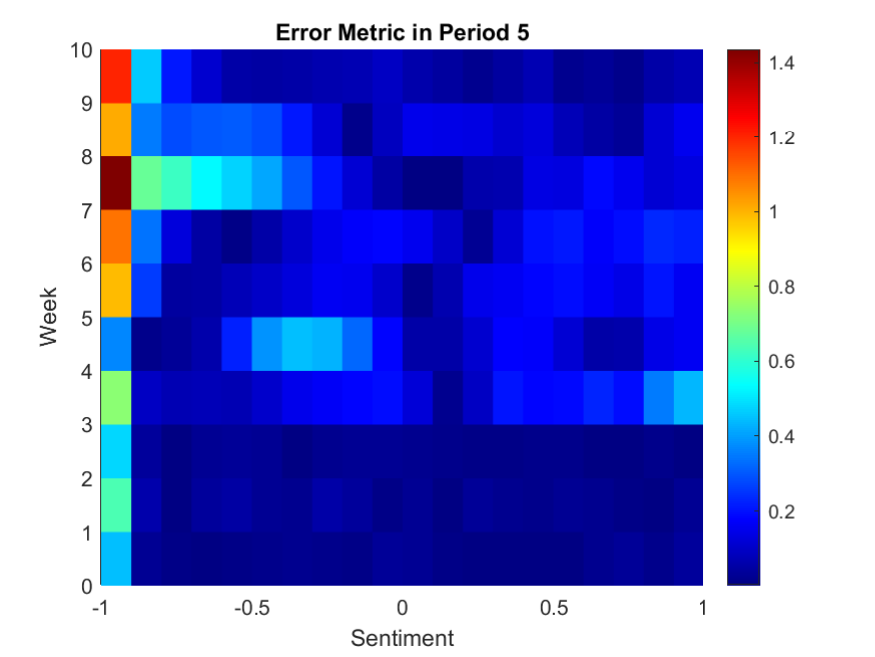} 
\includegraphics[page=6, width=0.325\textwidth]{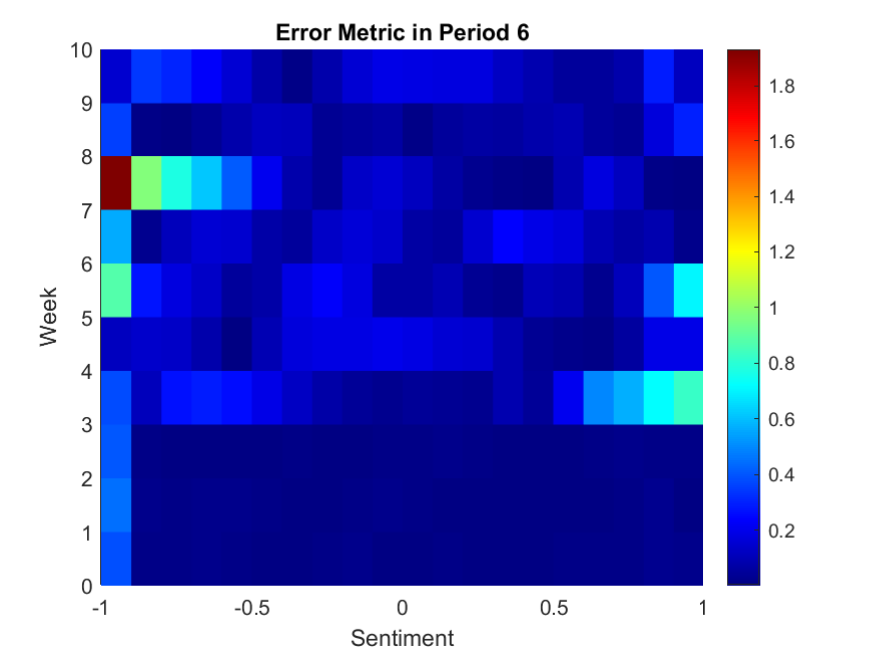} \\ \vspace{5mm}
\includegraphics[page=7, width=0.325\textwidth]{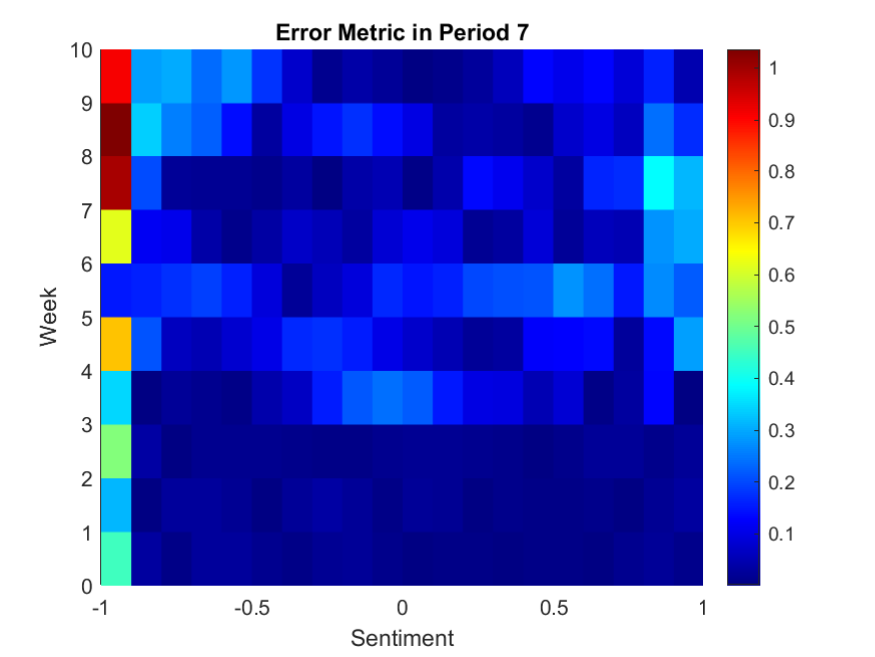}
\includegraphics[page=8, width=0.325\textwidth]{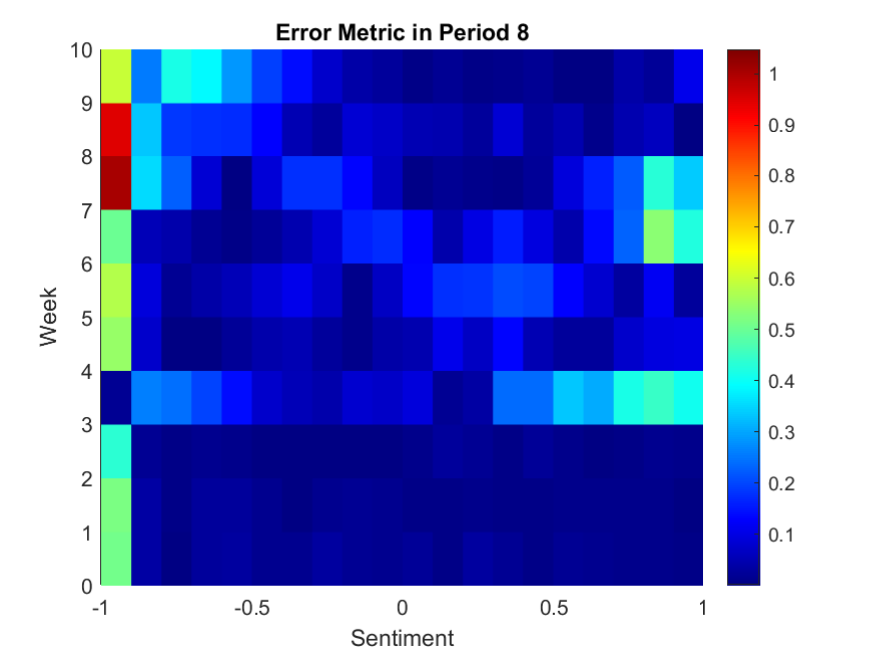} \\ \vspace{5mm}
\includegraphics[page=9, width=0.325\textwidth]{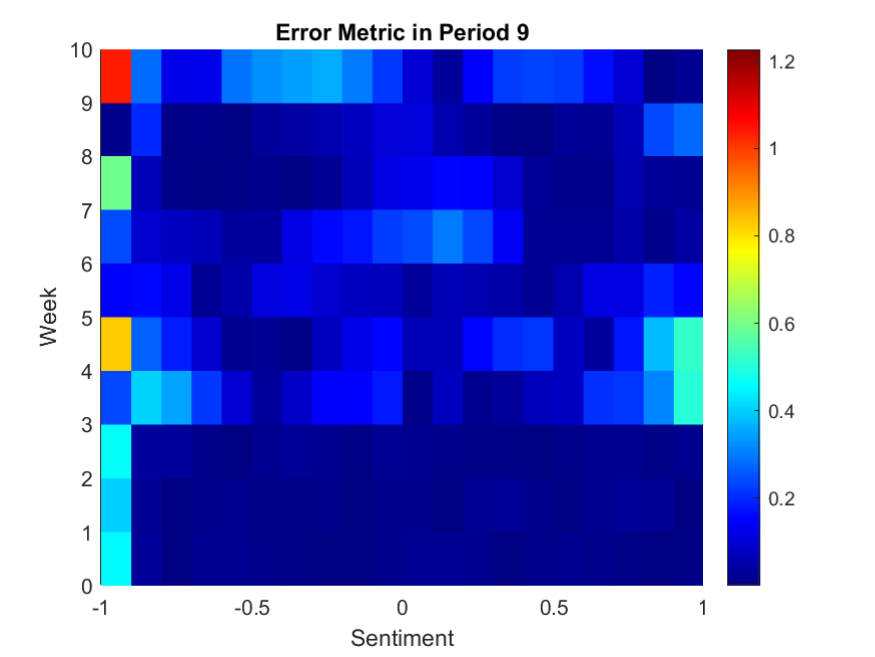} 
\includegraphics[page=10, width=0.325\textwidth]{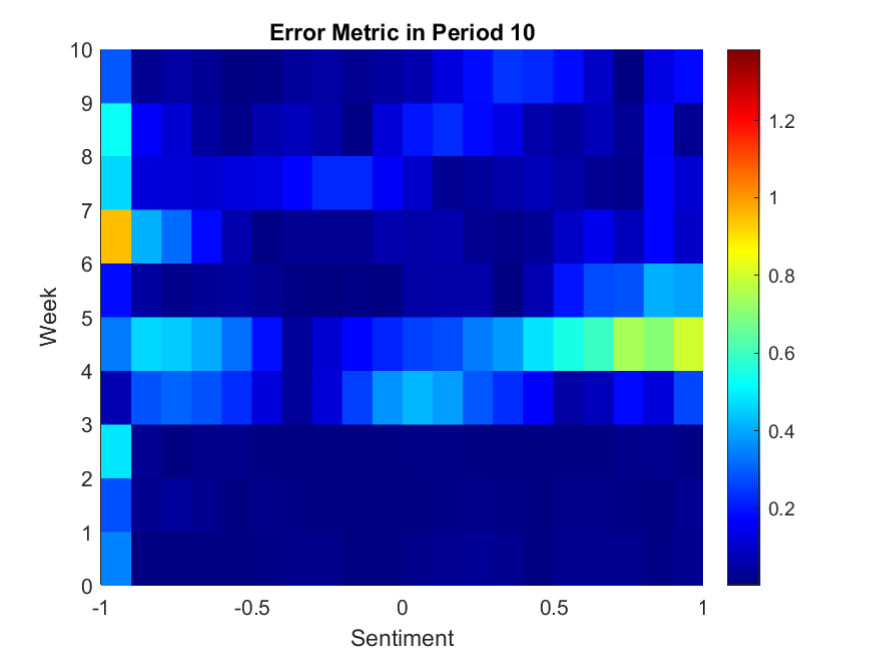}
\caption{Error metric \eqref{error_metric} for all ten forecasting periods.}
\label{fig:errors}
\end{figure}






\bibliographystyle{elsarticle-num-names} 
\bibliography{references}

@article{Achdou,
  author = {Y. Achdou and F. J. Buera and J-M. Lasry and P. L. Lions and B. Moll},
  title = {Partial differential equation models in macroeconomics},
  journal = {Philosophical Transactions of the Royal Society A: Mathematical, Physical and Engineering Sciences},
  volume = {372},
  number = {2028},
  pages = {20130397},
  year = {2014}
}

@article{Achdou2,
  author = {Y. Achdou and J. Han and J-M Lasry and P-L Lions and B. Moll},
  title = {Income and Wealth Distribution in Macroeconomics: A Continuous-Time Approach},
  journal = {The Review of Economic Studies},
  volume = {89},
  number = {1},
  pages = {45--86},
  month = {January},
  year = {2022}
}

@inproceedings{Banez,
  author = {R. A. Banez and H. Gao and L. Li and C. Yang and Z. Han and H. V. Poor},
  title = {Belief and opinion evolution in social networks based on a multi-population mean field game approach},
  booktitle = {2020 IEEE International Conference on Communications},
  pages = {1--6},
  year = {2020},
  doi = {10.1109/ICC40277.2020.9148985}
}

@article{Bauso,
  author = {D. Bauso and T. Hamidou and T. Basar},
  title = {Opinion dynamics in social networks through mean-field games},
  journal = {SIAM Journal on Control and Optimization},
  volume = {54},
  number = {6},
  pages = {3225--3257},
  year = {2016}
}

@article{BukhKlib,
  author = {A. L. Bukhgeim and M. V. Klibanov},
  title = {Uniqueness in the large of a class of multidimensional inverse problems},
  journal = {Soviet Math. Doklady},
  volume = {17},
  pages = {244--247},
  year = {1981}
}

@article{Cardaliaguet,
  author = {P. Cardaliaguet and C-A Lehalle},
  title = {Mean field game of controls and an application to trade crowding},
  journal = {Mathematics and Financial Economics},
  volume = {12},
  number = {3},
  pages = {335--363},
  year = {2018}
}

@article{Dogbe,
  author = {C. Dogb\'{e}},
  title = {Modeling crowd dynamics by the mean-field limit approach},
  journal = {Mathematical and Computer Modelling},
  volume = {52},
  number = {9-10},
  pages = {1506--1520},
  year = {2010}
}

@article{Festa,
  author = {A. Festa and S. G\"{o}ttlich and M. Ricciardi},
  title = {Forward-forward mean field games in mathematical modeling with application to opinion formation and voting models},
  journal = {Dyn. Games Appl.},
  volume = {15},
  pages = {2025},
  year = {2022},
  doi = {10.1007/s13235-024-00578-3}
}

@inproceedings{Firoozi,
  author = {D. Firoozi and P. E. Caines},
  title = {The execution problem in finance with major and minor traders: A mean field game formulation},
  booktitle = {International Symposium on Dynamic Games and Applications},
  publisher = {Springer International Publishing},
  address = {Cham},
  pages = {107--130},
  year = {2016},
  doi = {10.1007/978-3-319-70619-1_5}
}

@article{Gao,
  author = {H. Gao and A. Lin and R. A. Banez and W. Li and Z. Han and S. Osher and H. V. Poor},
  title = {Opinion evolution in social networks: Connecting mean field games to generative adversarial nets},
  journal = {IEEE Trans. Network Sci. Eng.},
  volume = {9},
  pages = {2734--2746},
  year = {2022},
  doi = {10.1109/TNSE.2022.3169057}
}

@article{HuangCainesMalhame2006,
  author = {M. Huang and R. P. Malham\'{e} and P. E. Caines},
  title = {Large population stochastic dynamic games: closed-loop McKean-Vlasov systems and the Nash certainty equivalence principle},
  journal = {Communications in Information and Systems},
  volume = {6},
  number = {3},
  pages = {221--252},
  year = {2006}
}

@inproceedings{HuttoGilbert2014,
  author = {C. J. Hutto and E. Gilbert},
  title = {VADER: A parsimonious rule-based model for sentiment analysis of social media text},
  booktitle = {Proceedings of the International AAAI Conference on Web and Social Media},
  volume = {8},
  pages = {216--225},
  month = {May},
  year = {2014},
  doi = {10.1609/icwsm.v8i1.14550}
}

@article{Klib95,
  author = {M. V. Klibanov and O. V. Ioussoupova},
  title = {Uniform strict convexity of a cost functional for three dimensional inverse scattering problem},
  journal = {SIAM J. Mathematical Analysis},
  volume = {26},
  pages = {147--179},
  year = {1995}
}

@article{Klib97,
  author = {M. V. Klibanov},
  title = {Global convexity in a three-dimensional inverse acoustic problem},
  journal = {SIAM J. Mathematical Analysis},
  volume = {28},
  pages = {1371--1388},
  year = {1997}
}

@book{KL,
  author = {M. V. Klibanov and J. Li},
  title = {Inverse Problems and Carleman Estimates: Global Uniqueness, Global Convergence and Experimental Data},
  publisher = {De Gruyter},
  address = {Berlin},
  year = {2021},
  doi = {10.1515/9783110745481}
}

@article{MFGSAM,
  author = {M. V. Klibanov and J. Li and H. Liu},
  title = {H\"{o}lder stability and uniqueness for the mean field games system via carleman estimates},
  journal = {Studies in Applied Mathematics},
  volume = {151},
  pages = {1447--1470},
  year = {2023},
  doi = {10.1111/sapm.12633}
}

@article{MFG1,
  author = {M. V. Klibanov and Y. Averboukh},
  title = {Lipschitz stability estimate and uniqueness in the retrospective analysis for the mean field games system via two Carleman estimates},
  journal = {SIAM J. Mathematical Analysis},
  volume = {56},
  pages = {616--636},
  year = {2024},
  doi = {10.1137/23M1554801}
}

@article{MFGCAMWA,
  author = {M. V. Klibanov and J. Li and Z. Yang},
  title = {Convexification for a coefficient inverse problem for a system of two coupled nonlinear parabolic equations},
  journal = {Computers and Mathematics with Applications},
  volume = {179},
  pages = {41--58},
  year = {2025},
  doi = {10.1016/j.camwa.2024.12.004}
}

@article{Kepid,
  author = {M. V. Klibanov and J. Li and Z. Yang},
  title = {Spatiotemporal monitoring of epidemics via solution of a coefficient inverse problem},
  journal = {Inverse Problems and Imaging},
  volume = {19},
  pages = {1142--1166},
  year = {2025}
}

@misc{Kliba2025,
  author = {M. V. Klibanov and K. McGoff and T. Truong},
  title = {Forecasting public sentiments via mean field games},
  note = {arXiv preprint; accepted in SIAM J. Applied Mathematics},
  year = {2025},
  url = {https://arxiv.org/abs/2506.08465}
}

@article{MFG8IPI,
  author = {M. V. Klibanov and J. Li and Z. Yang},
  title = {Convexification numerical method for a coefficient inverse problem for the system of nonlinear parabolic equations governing mean field games},
  journal = {Inverse Problems and Imaging},
  volume = {19},
  number = {2},
  pages = {219--252},
  year = {2025},
  doi = {10.3934/ipi.2024031}
}

@article{MFG2,
  author = {M. V. Klibanov},
  title = {The mean-field games system: Carleman estimates, Lipschitz stability and uniqueness},
  journal = {J. Inverse Ill-Posed Probl.},
  volume = {31},
  pages = {455--466},
  year = {2023},
  doi = {10.1515/jiip-2023-0023}
}

@article{Lachapelle,
  author = {A. Lachapelle and M-T Wolfram},
  title = {On a mean field game approach modeling congestion and aversion in pedestrian crowds},
  journal = {Transportation research part B: methodological},
  volume = {45},
  number = {10},
  pages = {1572--1589},
  year = {2011}
}

@article{LasryLions2007,
  author = {J.-M. Lasry and P.-L. Lions},
  title = {Mean field games},
  journal = {Japanese Journal of Mathematics},
  volume = {2},
  number = {1},
  pages = {229--260},
  year = {2007}
}

@book{MFGbook,
  author = {M. V. Klibanov and J. Li},
  title = {Carleman Estimates in Mean Field Games},
  publisher = {De Gruyter},
  year = {2025}
}

@inproceedings{Stella,
  author = {L. Stella and F. Bagagiolo and D. Bauso and G. Como},
  title = {Opinion dynamics and stubbornness through mean-field games},
  booktitle = {52nd IEEE Conference on Decision and Control},
  pages = {2519--2524},
  year = {2013},
  doi = {10.1109/CDC.2013.6760259}
}

@book{Tikhonov1995,
  author = {A. N. Tikhonov and A. V. Goncharsky and V. V. Stepanov and A. G. Yagola},
  title = {Numerical Methods for the Solution of Ill-Posed Problems},
  series = {Mathematics and Its Applications},
  publisher = {Springer},
  address = {Dordrecht},
  year = {1995},
  isbn = {978-0-7923-3583-2},
  doi = {10.1007/978-94-015-8480-7},
  note = {Originally published in Russian}
}

@misc{Matlab2024a,
year = {2024},
author = {The MathWorks Inc.},
title = {Optimization toolbox},
publisher = {The MathWorks Inc.},
address = {Natick, Massachusetts, United States},
url = {https://www.mathworks.com/help/optim/}
}

\end{document}